\documentclass[times]{MyArticle}

\usepackage[colorlinks,bookmarksopen,bookmarksnumbered,citecolor=red,urlcolor=red]{hyperref}

\usepackage[usenames,dvipsnames,svgnames,table]{xcolor}

\newcommand\BibTeX{{\rmfamily B\kern-.05em \textsc{i\kern-.025em b}\kern-.08em
T\kern-.1667em\lower.7ex\hbox{E}\kern-.125emX}}

\begin{document}

\title{A Meshfree Lagrangian Method for Flow on Manifolds}

\author{Pratik Suchde \affil{1}\corrauth}

\address{\affilnum{1}Fraunhofer ITWM, 67663 Kaiserslautern, Germany}

\corraddr{E-mail: pratik.suchde@itwm.fraunhofer.de}

\begin{abstract}
In this paper, we present a novel meshfree framework for fluid flow simulations on arbitrarily curved surfaces. First, we introduce a new meshfree Lagrangian framework to model flow on surfaces. Meshfree points or particles, which are used to discretize the domain, move in a Lagrangian sense along the given surface. This is done without discretizing the bulk around the surface, without parametrizing the surface, and without a background mesh. A key novelty that is introduced is the handling of flow with evolving free boundaries on a curved surface. The use of this framework to model flow on moving and deforming surfaces is also introduced. Then, we present the application of this framework to solve fluid flow problems defined on surfaces numerically. In combination with a meshfree Generalized Finite Difference Method (GFDM), we introduce a strong form meshfree collocation scheme to solve the Navier--Stokes equations posed on manifolds. Benchmark examples are proposed to validate the Lagrangian framework and the surface Navier--Stokes equations with the presence of free boundaries. 
\end{abstract}


\maketitle

\vspace{-6pt}


\section{Introduction}
\label{sec:Introduction}

Fluid flow on curved surfaces forms an important aspect in modelling dynamics on biomembranes and other interfacial transport processes \cite{Barrett2015, Edwards1991}, in representing surfactant dynamics \cite{Grassia2016}, in thin film flow and lubrication \cite{OBrien2002}, and visualization \cite{Stam2003}, among other applications. 

One possibility for simulating flow on curved surfaces is the Navier--Stokes equations posed on manifolds. Several different and unequivalent formulations of the Navier--Stokes equations on manifolds have appeared in literature. This is mainly due to the different interpretations of the viscous stress tensor. This, in turn, arises from the different possibilities of the Laplace operator on vector fields on manifolds \cite{Chan2017}. Using different arguments, recent work by \cite{Chan2017, Jankuhn2017, Koba2017, Miura2017} clears this confusion to show that the physically accurate interpretation is one based on the so-called Boussinesq-Scriven surface stress tensor \cite{Boussinesq1913, Scriven1960}, which has also been referred to as the deformation tensor \cite{Ebin1970}. For a comprehensive review of the different interpretations of the Navier--Stokes equations on manifolds, we refer to \cite{Chan2017} for Riemannian manifolds, and \cite{Jankuhn2017} for the case of the surface Navier--Stokes equations, which we consider in the present work.

Most numerical methods to solve surface flow have used simplified models. Visualization driven work has  mostly considered only inviscid flow \cite{Auer2012, Auer2013}. A lot of work has been done for surface flow on surfaces of specific shapes, such as spheres \cite{Fengler2005} or other radial manifolds \cite{Gross2018}. Atmospheric dynamics work to model flow on the surface of the Earth solve two-dimensional simulations in the latitude-longitude frame \cite{Dritschel2015, Qi2014}. Generalizing this to arbitrary geometries would involve parametrizing the surface, which could prove to be very challenging for complex or time-varying geometries, or for flow involving free boundaries. 

Numerical methods to solve the full incompressible surface Navier--Stokes equations has gained an increasing interest in the last few years \cite{Fries2018, Olshanskii2019, Reuther2019}. Of particular interest in the present context is the recent work that recasts surface flow into a tangential differential calculus framework \cite{Fries2018, Jankuhn2017}. The advantage of this framework is that the equations are formulated entirely in Eucledian space, including the derivative computations, and thus they are easily implemented. As a result, it avoids the requirement of maintaining a parametrization of the surface, and the related issues with singularities arising in metrics. This framework is generally used to solve the underlying PDEs in a weak formulation. We adopt this framework, and use it to solve the PDEs in a strong form. Here, we focus on the use of this framework in an intrinsic setting, where the underlying PDEs are solved directly on the surface. Under the terminology of `calculus on surfaces without parametrization', similar methods have also been used in embedding methods which solve PDEs on a band around the surface, such as the closest point method \cite{Marz2012,Ruuth2008}.


Solutions to the surface Navier--Stokes equations typically rely on Finite Element methods \cite{Fries2018, Nitschke2012, Olshanskii2019, Reuther2019}. Here, we introduce the use of a meshfree method to solve this problem. In the last two decades, meshfree methods have started to become a popular alternative for conventional flow applications on volume domains, especially when the computational domain is complex or time-dependent. In such situations, mesh generation is very complex, and often needs several man hours. Furthermore, rapidly evolving domains result in the need of extensive remeshing which can be a very expensive portion of the simulation time. On the other hand, point cloud generation for meshfree methods has been largely automated, and is thus much faster for complex domains. Additionally, fixing distortion in point clouds can be done locally, and as a result, the meshfree equivalent of remeshing is also much cheaper. Many of these advantages carry over to the case of surface flow. The extensive need for expensive volume mesh generation for complex geometries is significantly reduced for the case of surface meshes. However, the issue of costly fixing of mesh distortion for time-dependent geometries is still present. This is relevant when the surface itself is evolving, and when there are evolving free boundaries within a surface. In our earlier work \cite{Suchde2019_MovingSurfaces}, we have shown the advantage of meshfree Lagrangian frameworks for the former case to solve PDEs on evolving surfaces. Here, we explore their use in the latter in the context of flow on surfaces. 

We introduce a novel meshfree Lagrangian framework to model flow with free boundaries on a moving curved surface. Meshfree points undergo Lagrangian motion with two velocities: the velocity of the fluid, and that of the surface. A point cloud is used to represent not only the moving fluid, but also the surface constraining the flow. Another novelty in the present work is the inclusion of free boundaries in the context of Navier--Stokes on manifolds, which has not been done before to the best of our knowledge.


The paper is organized as follows. Section~\ref{sec:Prelim} introduces some basic notions of surface points clouds. Section~\ref{sec:Lagrangian} presents a novel Lagrangian framework to represent moving fluid points on a moving surface. This is then extended to moving free boundaries on a possibly moving surface. Section~\ref{sec:DiffOperators} provides an overview of surface differential operators and numerical methods to discretize them in a Generalized Finite Difference Method~(GFDM) setting. The Navier--Stokes equations on manifolds are introduced in Section~\ref{sec:NSman}, along with a numerical scheme to solve them using a strong-form GFDM method. Section~\ref{sec:Results} then presents some numerical results for the surface Navier--Stokes equations, and the paper is concluded with a brief discussion in Section~\ref{sec:conclusion}. 


\section{Preliminaries}
\label{sec:Prelim}

We consider a smooth orientable surface $M \subset \mathbb{R}^3$, which may or may not have boundaries. The notations used for surface point clouds in this paper follow \cite{Suchde2019_MovingSurfaces}, and those for the surface flow problem follow \cite{Fries2018}. We start with the case of point cloud surfaces without free boundaries. The inclusion of free boundaries is done in Section~\ref{sec:FreeBoundary}. 

\begin{figure}
  \centering
  \includegraphics[width=0.7\textwidth,trim={1cm 2cm 8cm 5cm},clip]{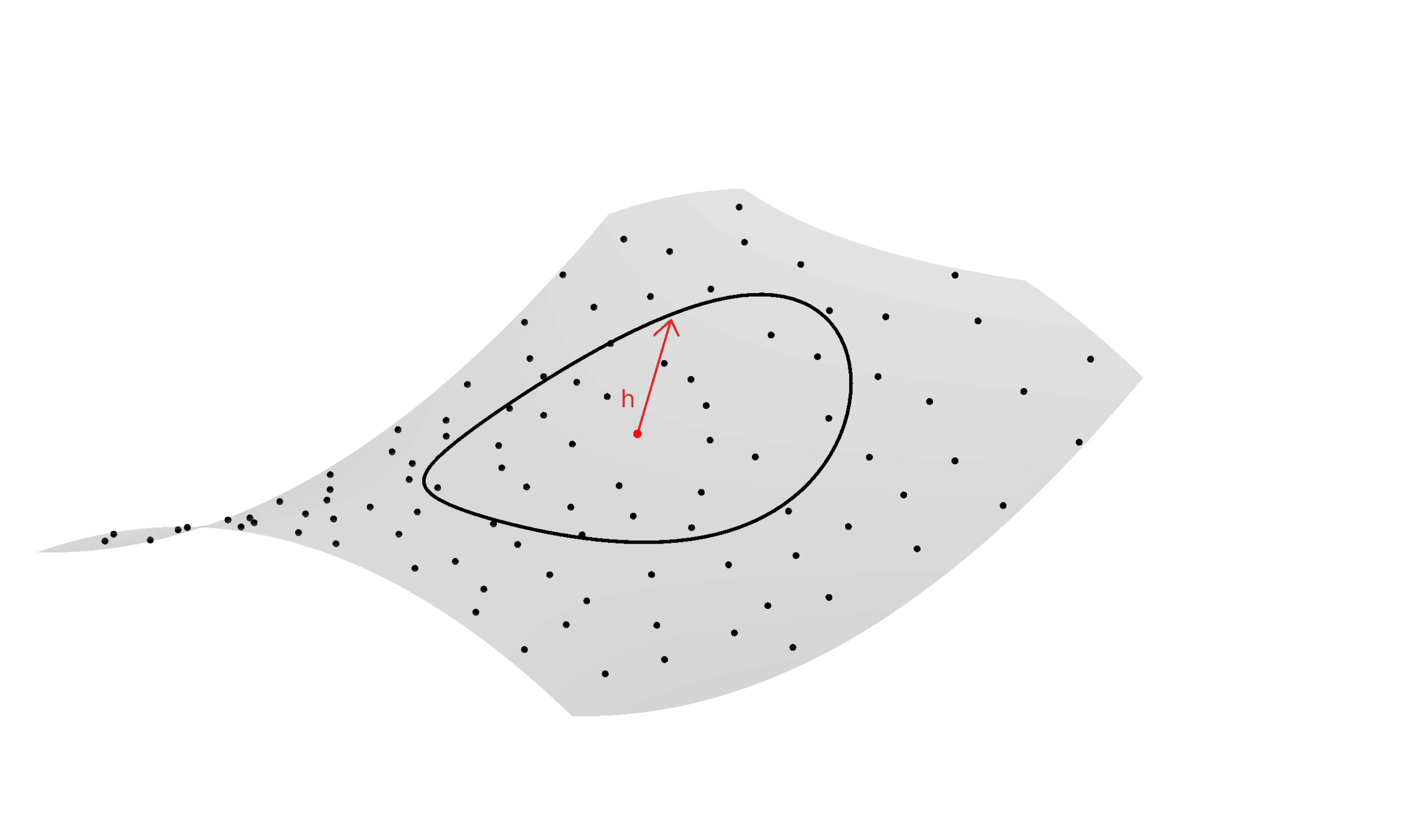} 
  \caption{Surface neighbourhood determined by distances in the embedding space, $\mathbb{R}^3$. The `circle' shown demarcating the neighbourhood is the intersection of the sphere of radius $h$ with the surface.}
  \label{Fig:SurfaceNeighbourhood}%
\end{figure}

The surface is discretized with a point cloud $PC$ consisting of $N = N(t)$ irregularly spaced points. Approximations at each point are based on a compact support or neighbourhood consisting of all points within a distance $h = h(\vec{x}, t)$ from it. An example of such a neighbourhood is shown in Figure~\ref{Fig:SurfaceNeighbourhood}. All distances are computed as standard Euclidean distances, and not distances along the surface. $h$ is referred to as the smoothing length or interaction radius. Efficient neighbour searching algorithms for volume-domain meshfree methods (for example, \cite{Dominguez2010, Drumm2008, Onderik2008}) directly carry over to the surface case here.

Inter-particle distance ranges from $r_{min}h$ to $r_{max}h$, for fixed simulation-independent parameters $r_{min}$ and $r_{max}$. The parameters used here are adopted from conventions in volumetric meshfree flow simulations \cite{Drumm2008, Jefferies2015, Suchde2017_CCC}. We set $r_{min} = 0.2$ and $r_{max} = 0.45$, as has been done for meshfree surfaces in \cite{Suchde2019_StaticSurfaces, Suchde2019_MovingSurfaces}. This results in about $15-20$ points in each neighbourhood. This set up leads to $h$ being a measure of not just the support size, but also the discretization size. On a moving point cloud, these minimum and maximum inter-particle distances are ensured by inserting and removing points. More details of this are explained in Section~\ref{sec:Distortion}.

\begin{figure}
  \centering
  \includegraphics[width=0.32\textwidth]{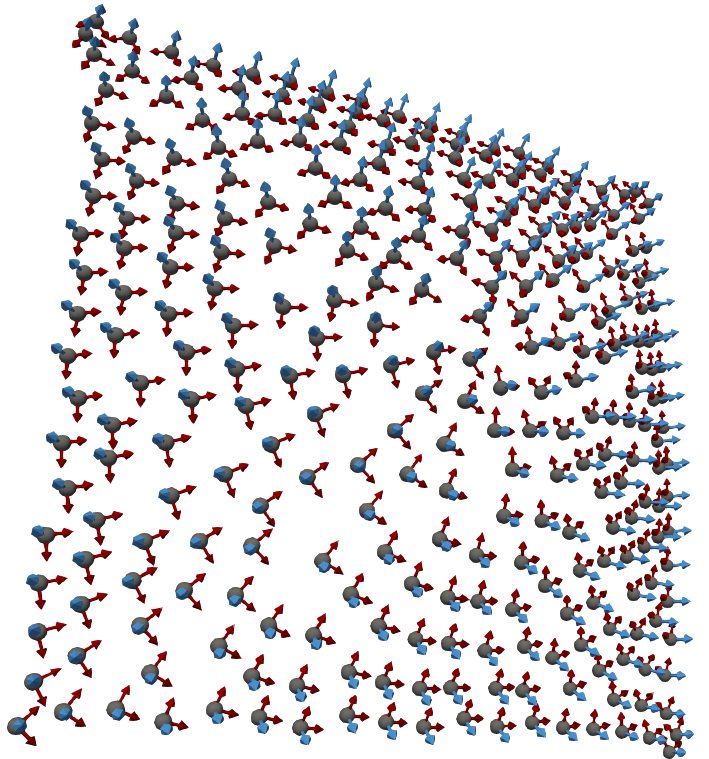}
  \includegraphics[width=0.32\textwidth]{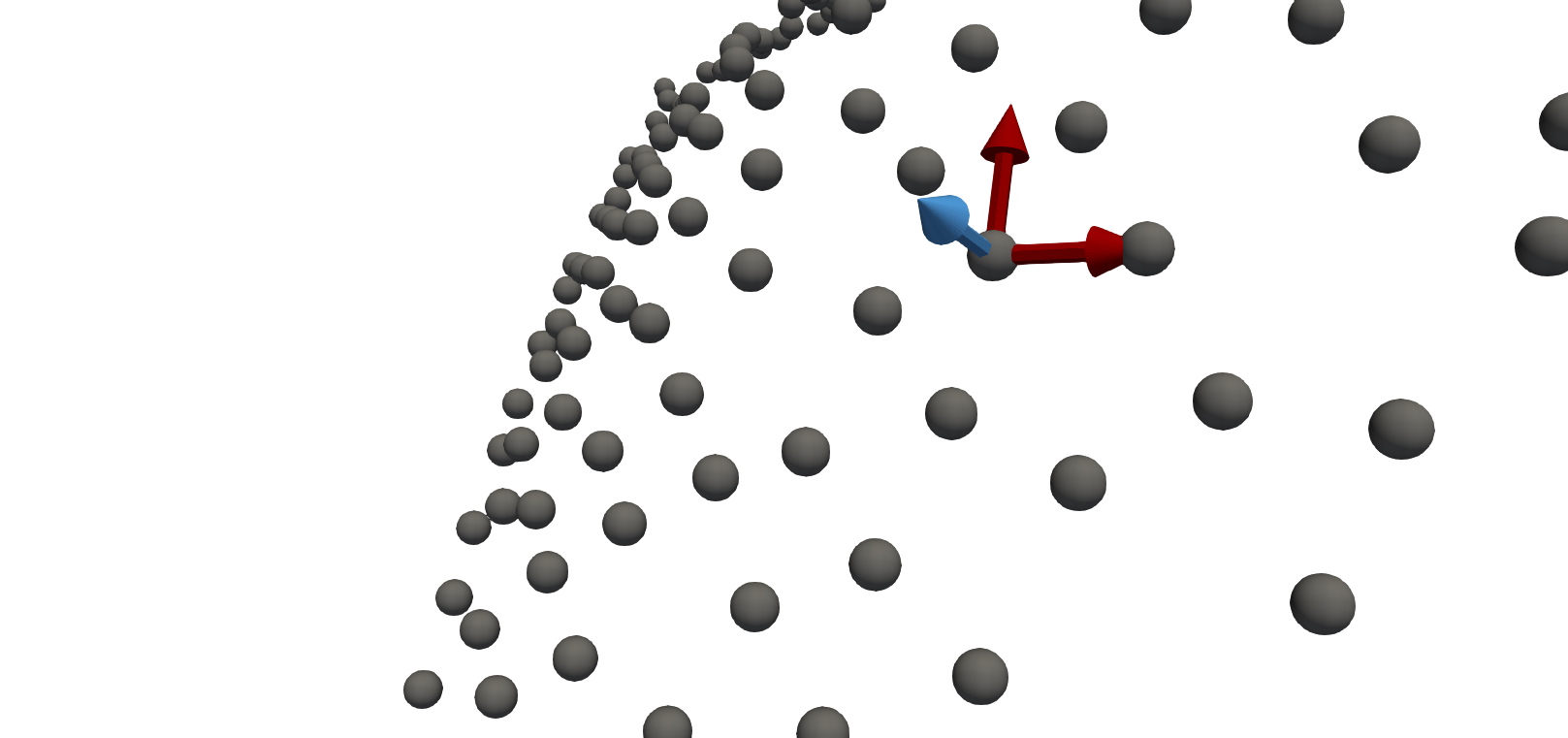}
  \includegraphics[width=0.32\textwidth]{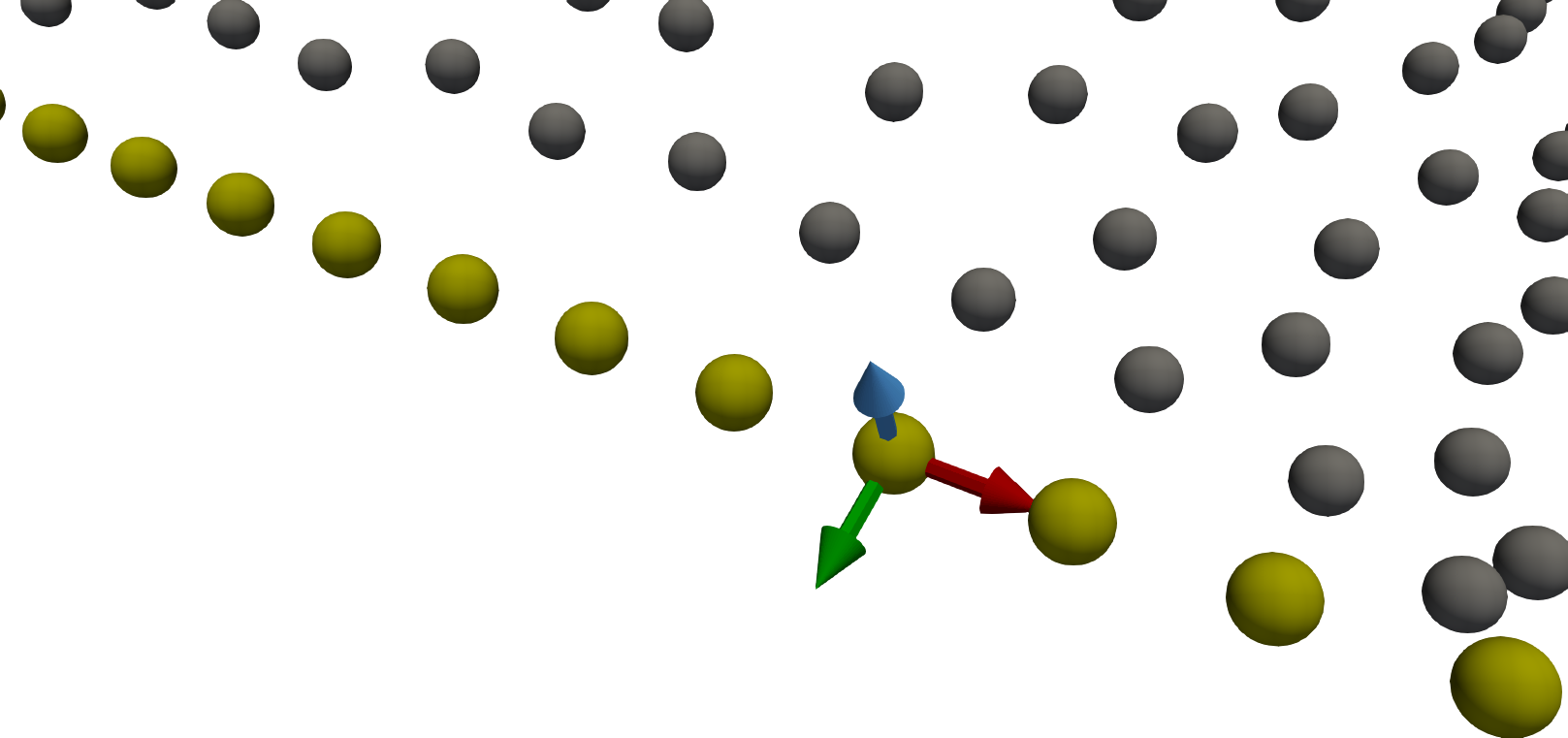}  
  \caption{Normal and tangents to the surface. The center figure shows the case of an interior point. The right figure shows the same for a boundary point. The surface normal is shown in blue, the boundary normal is shown in green, and the tangent(s) in red.}
  \label{Fig:NormalTangent}%
\end{figure}

Each point is equipped with a unit surface normal $\vec{n}$. It is ensured that the normals are defined such that the discrete surface is oriented. Boundary points further have an outward pointing unit boundary normal $\vec{\nu}$, which is conormal to $\vec{n}$. This is illustrated in Figure~\ref{Fig:NormalTangent}. Details of computation of the normals on moving point clouds follow the procedure lined out in \cite{Suchde2019_MovingSurfaces}. 
%


\section{The Lagrangian Framework}
\label{sec:Lagrangian}

In our earlier work \cite{Suchde2019_MovingSurfaces}, we introduced a meshfree Lagrangian framework to model evolving-in-time surfaces numerically, including meshfree contact handling algorithms. Here, we extend that framework to model fluid flow within a (possibly moving) surface. We propose to do this using Lagrangian particles which move along the given surface. This introduces several special features, as we explain below.

We emphasize that in the Lagrangian framework that follows, we do not maintain any correlation to the domain at the initial state. At a specific time step, all approximations are performed on the point cloud at that time step directly, without mapping it to the initial configuration. 

For the case of fluid motion on a moving surface, two advection velocities are relevant. The velocity $\vec{v}$ of the fluid moving on the surface, and the velocity $\vec{w}$ of the surface itself. Thus, a particle on the surface moves with a velocity of $\vec{v}+\vec{w}$, and the movement of the point cloud is given by 
\begin{equation}
	\label{sec:Move}
	\frac{d\vec{x}}{dt} = \vec{v} + \vec{w}\,.
\end{equation}
We note that $\vec{v}$ is necessarily tangential to the surface, i.e., $\vec{v}\cdot \vec{n}=0$. Whereas $\vec{w}$ is not constrained and can have components both normal and tangential to the surface. In \cite{Suchde2019_MovingSurfaces}, we only handled the special case of $\vec{v}=0$ and $\vec{w}\neq 0$, which was used to solve PDEs on evolving surfaces, in the absence of any flow on the surface. Here, we deal with the general case of $\vec{v}\neq 0$  and $\vec{w} \neq 0$.

We start by restricting ourselves to flow without free boundaries, handling of which is introduced in Section~\ref{sec:FreeBoundary}.

\subsection{The Lagrangian movement of points}
\label{sec:LagMove}

On the numerical level, we split the movement into two steps. The first is based on the fluid velocity, and the second one is based on the surface velocity. Each step involves a second order in time movement. We note that most Lagrangian methods only use a first order movement, which has been shown to be extremely inaccurate in capturing rotational components of motion \cite{Suchde2018_PCM}. For time integration between time levels $t^n$ and $t^{n+1}$ with the velocity $\vec{v}$, we get
\begin{equation}
	\label{Eq:MoveV}
		\Delta \vec{x}_1 = \vec{v}^{\,(n)}\Delta t + \frac{1}{2}\frac{\vec{v}^{\,(n)} - \vec{v}^{\,(n-1)}}{\Delta t_0} (\Delta t)^2 \,,
\end{equation}
where bracketed superscripts indicate the time level, $\Delta t = t^{n+1} - t^{n}$ is the current time step, and $\Delta t_0 = t^{n} - t^{n-1}$ is the previous time step. Here, we assume that $\vec{v}^{\,(n+1)}$ is unknown, as is the case if the movement step is done before computing the new velocity. If $\vec{v}^{\,(n+1)}$ is known during the movement step, $\vec{v}^{\,(n)}$ and $\vec{v}^{\,(n-1)}$ can be replaced by $\vec{v}^{\,(n+1)}$ and $\vec{v}^{\,(n)}$ respectively in Eq.\,\eqref{Eq:MoveV}.

Since the velocity $\vec{v}$ lies in the tangent space, an exact integration along the velocity streamlines would result in particles still being on the surface. However, the discrete movement based on Eq.\,\eqref{Eq:MoveV} can cause points to be moved off the surface to a distance of the order of $\mathcal{O}(\Delta t ^3)$. To avoid accumulation of these numerical errors in this movement step, we follow it with a projection of the point cloud back to the surface given by the point cloud at $t^n$, before the movement. Details of the projection step follow that in \cite{Suchde2019_MovingSurfaces}. 

Eq.\,\eqref{Eq:MoveV} and the projection step are then followed by movement according to the surface velocity $\vec{w}$
\begin{equation}
	\label{Eq:MoveW}
		\Delta \vec{x}_2 = \vec{w}^{\,(n)}\Delta t + \frac{1}{2}\frac{\vec{w}^{\,(n)} - \vec{w}^{\,(n-1)}}{\Delta t_0} (\Delta t)^2 \,.
\end{equation}
Thus, the point cloud at the new time level is given by
\begin{equation}
	\label{Eq:MoveOverall}
	\vec{x}^{\,(n+1)} = \tilde{P}_{(n)} \left( \vec{x}^{\,(n)} + \Delta \vec{x}_1 \right) + \Delta \vec{x}_2 \,,
\end{equation}
for the projection step $\tilde{P}_{(n)}$ to the point cloud at $t^n$. 

We perform a second order movement step with each velocity to more accurately capture the movement. However, it must be noted that the splitting explained above is only of first order. 

\subsection{Fixing Distortion}
\label{sec:Distortion}

The ideas of fixing distortion in surface point clouds directly carry over from \cite{Suchde2019_MovingSurfaces}. We present a brief description here for the sake of completeness. As introduced in Section~\ref{sec:Prelim}, to maintain regularity of the point cloud, it is ensured that no two points come closer than a distance of $r_{min}h$, and that there is no hole in the point cloud of radius $r_{max}h$ that contains no points. As the point cloud is moved in the Lagrangian sense, this regularity condition can be violated, resulting in a distorted point cloud. Detecting and fixing these violations is always done entirely locally, unlike the case of fixing a distorted mesh which could require a global remeshing. Violations to the closeness condition can be detected trivially with distance computations within each neighbourhood. If two points come closer than than $r_{min}h$ apart, they are merged into one point at the center location, with all physical properties interpolated at the new location. Detection of formed holes in the point cloud are done by checking the circumradius of triangles formed by local tessellations within a compact support. Since $r_{min}$ and $r_{max}$ are fixed parameters, changing $h$ automatically triggers the addition and merging algorithms, and can thus be used to introduce adaptive refinement. Detection of holes and close points is done immediately after the Lagrangian movement step. First, the distortion fixing is done specific to boundary points, after which it is done for interior points. We refer to Section 3 of \cite{Suchde2019_MovingSurfaces} for more details. Fixing distortion in this manner ensures that there is no need for any artificial movement step to maintain point cloud regularity, as is the case in many Lagrangian particle methods like SPH (for example, \cite{Pahar2016}). 


\subsection{Lagrangian Framework with Free Boundaries}
\label{sec:FreeBoundary}

We now extend the above notions of flow on point cloud surfaces to include free boundaries. In the above case, the fluid region being simulated occupies the entire computational domain. i.e. $PC$ describes not just the fluid, but also the surface that the fluid is restricted to. However, for the case of free boundaries, this is no longer true. There, the fluid occupies only a part of the surface domain. Thus, there is an additional need to represent the surface which constrains the fluid movement. 

\begin{figure}
	\centering
	\includegraphics[width=0.55\textwidth]{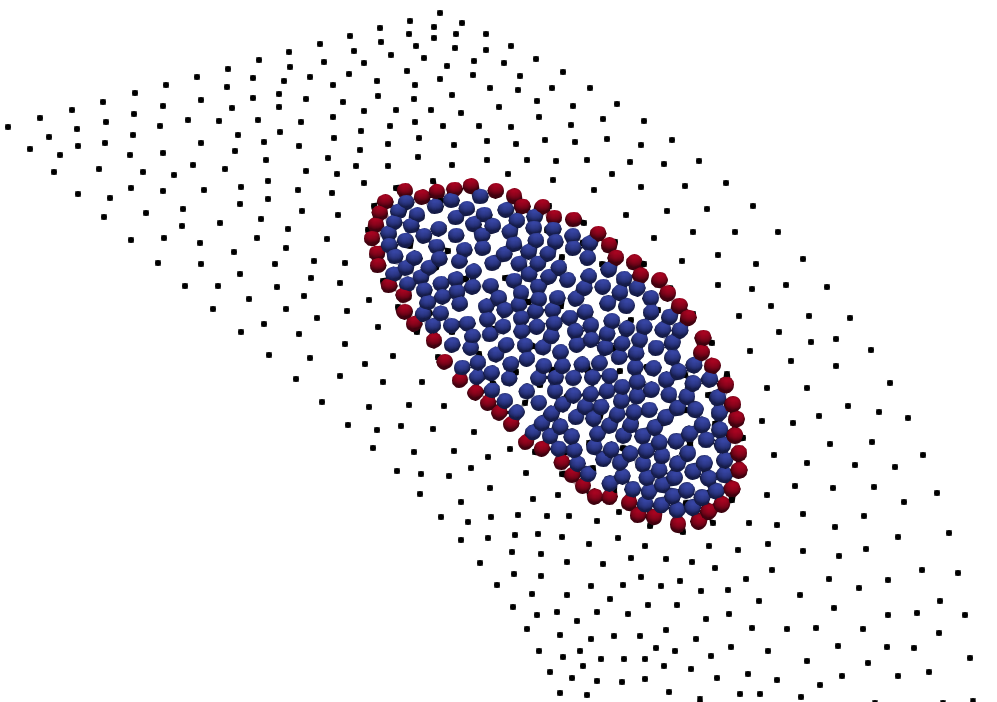} 
	\caption{For simulating free boundaries in surface flow, two point clouds are maintained. The main point cloud representing the fluid is shown in spherical glyphs: blue for interior points, and red for free boundary points. The background point cloud representing the shape of the overall computational domain is shown with black points. }
	\label{Fig:FreeBoundary}%
\end{figure}

This needed surface discretization could be done in multiple ways. One possible method is to maintain a surface mesh, over which the particles move. This would represent the extension of particle-in-cell methods \cite{Evans1957} to surface domains. This approach can be very useful if the movement of the surface is restricted to rigid body motion. However, for time-evolving surfaces undergoing arbitrary movements, the movement of the mesh could cause distortion which may require expensive remeshing to fix. Thus, the advantages of point clouds again become relevant, and we leverage them by using a point cloud even for the surface discretization.  To incorporate free boundaries, we introduce the use of two different point clouds. A main point cloud $PC$ representing the moving fluid, which moves over a background point cloud $PC^{\text{back}}$ representing the surface. The background point cloud is only used to constrain the movement of the main point cloud. Underlying PDEs are only solved on $PC$ and not on $PC^{\text{back}}$. The use of two points clouds in such a manner is illustrated in Figure~\ref{Fig:FreeBoundary}.

Notations for the main point cloud  $PC$ follow that introduced in Section~\ref{sec:Prelim} with $N = N(t)$ points, and a smoothing length of $h=h(\vec{x}, t)$. The background point cloud $PC^{\text{back}}$ to define the surface contains $N^{\text{back}} = N^{\text{back}}(t)$ points, and a smoothing length of $h^{\text{back}}=h^{\text{back}}(\vec{x}, t)$. 

$PC^{\text{back}}$ is moved only with the surface velocity $\vec{w}$, analogously to Eq.\,\eqref{Eq:MoveW}. If the surface is stationary ($\vec{w}=0$), the background point cloud does not move, and thus does not change in time. The main point cloud is moved with both $\vec{v}$ and $\vec{w}$, as described in Section~\ref{sec:LagMove}, with the only difference being the projection. In presence of the background point cloud, the projection after the move with the fluid velocity is done to the surface given by $PC^{\text{back}}$. Thus the movement can be summarized as follows. First, the locations of the main fluid point cloud is updated according to Eq.\,\eqref{Eq:MoveV}. This is then projected to the surface given by the background point cloud. Then both point clouds are moved with the surface velocity according to Eq.\,\eqref{Eq:MoveW}. 

Distortion needs to be fixed for both point clouds upon movement. For both, this follows the same procedure outlined in Section~\ref{sec:Distortion}. 

For ease of visualization, in some of the figures to follow, we represent the background point cloud with a tessellation of it. Henceforth, unless mentioned otherwise, the term point cloud will refer to the main point cloud representing the fluid. We emphasize here that the background point cloud is only maintained for applications involving free boundaries. Without free boundaries within the surface, there is no need to maintain $PC^{\text{back}}$, as the main point cloud sufficiently describes both the fluid and the underlying surface.  

\subsection{Free Boundary Detection}

An important step in simulating a fluid with free boundaries is the detection of the points that represent the free boundary. Initially, the prescription of free boundary points is straight forward and can be done directly based on the initial fluid configuration. However, extra care needs to be taken as the free boundary evolves, to handle events such as the fluid domain hitting a fixed boundary, where the distinction between free and fixed boundary points needs to be made for the correct application of the boundary conditions. An important use case is when a topological change occurs within the fluid domain, for example when a fluid droplet splits or merges with another. 

\begin{figure}
  \centering
  \includegraphics[trim={10cm 2cm 8cm 5cm},clip,width=0.49\textwidth]{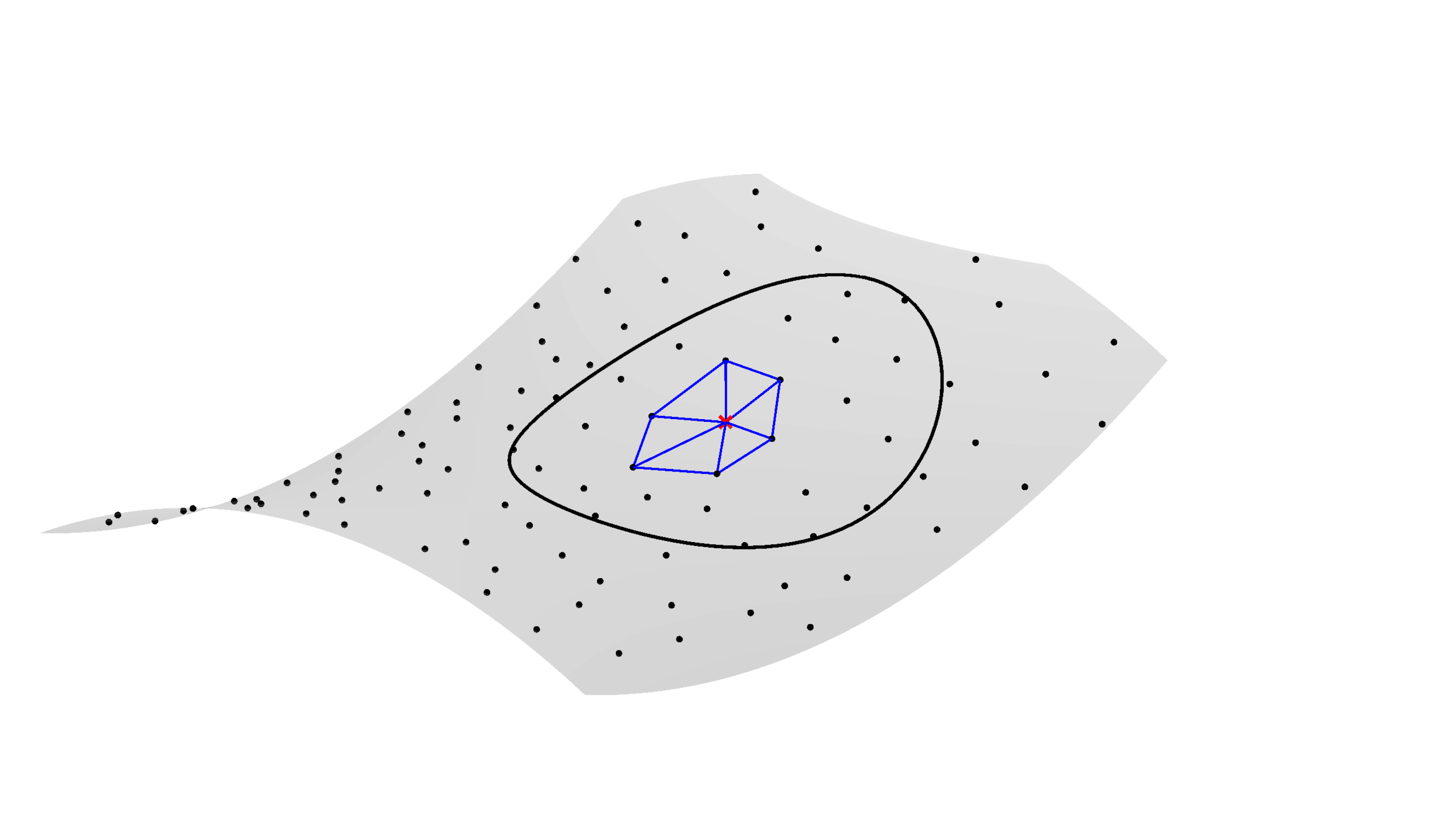}   
  \includegraphics[trim={8cm 2cm 15cm 5cm},clip,width=0.49\textwidth]{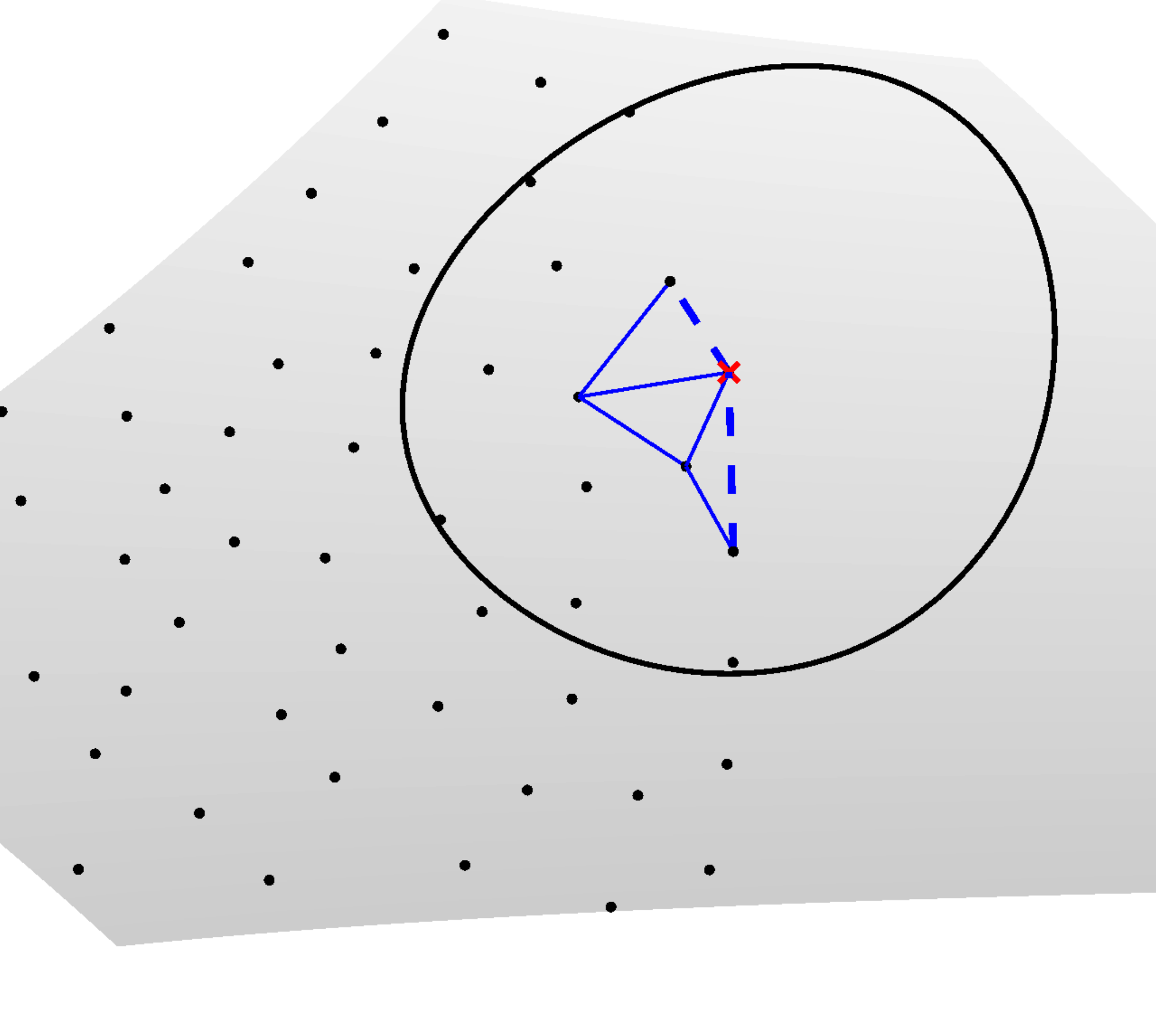}
  \caption{Local tessellation used for free boundary detection and fixing distortion. Interior point~(left). Boundary point~(right) with the open edges marked with thicker dashed lines.}
  \label{Fig:LocalTri}%
\end{figure}

Several approaches have been proposed for detecting free surfaces in volumetric flow (e.g.,  \cite{Marrone2010,Saucedo2019, Tiwari2003}) , and those can be extended to the present context. We choose to reuse the local tessellations which are already being used to fix distortion. For each point, a ``1-ring" of triangles is determined by tessellating its neighbourhood, as shown in Figure~\ref{Fig:LocalTri}. We emphasize here that these ``1-rings" of triangles at every point need not stitch together to form a global mesh on the point cloud \cite{Suchde2017_CCC}. Furthermore, no restriction on the quality of these triangles is enforced, nor are any approximations done based on these triangles. Once the tessellation is computed, a point is labelled to be part of the free boundary if there are open edges incident on the point. An example of this is illustrated in Figure~\ref{Fig:LocalTri}. There is an additional possibility of a fixed boundary point for points that lie on boundaries of the surface.

In computing the local tessellations, an important heuristic that is followed is to avoid forming a triangle between three boundary points of the previous time step. For boundary point detection, this is needed to prevent free boundary points from being incorrectly labelled as interior points. Since this tessellation is also used to find holes in the point cloud, avoiding triangles between boundary points of the previous time step ensures that regions outside the fluid domain are not incorrectly identified as holes, and no points are added outside the actual domain of computation. A test case which highlights the use of these algorithms is presented later, in Section~\ref{sec:DropFall}.

\subsection{Contact handling}
The contact handling algorithms between meshfree surface point clouds \cite{Suchde2019_MovingSurfaces} directly carry over to the present case. These are relevant when free boundaries are present, especially when the topology of the fluid domain changes. For example, when a fluid region splits into multiple disconnected regions, or when movement of a free boundary causes a ``hole" to be formed in the middle of the fluid. 

\subsection{Numerical Examples}
\label{sec:LagNum}

To verify the Lagrangian framework introduced, we consider a few examples of a fluid with free boundaries moving on an evolving surface, with prescribed velocities $\vec{v}$ and $\vec{w}$. We emphasize that only the ODE for the Lagrangian movement of the point cloud is solved here. 

\subsubsection{Moving droplet on a fixed cylinder}
\begin{figure}
	\centering
	\includegraphics[width=0.4\textwidth]{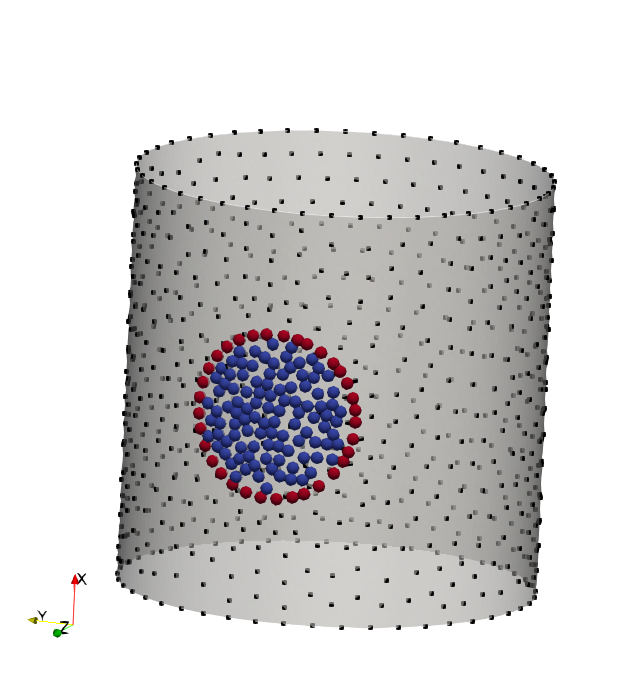} 
	\caption{Initial condition for testing the Lagrangian motion of a free fluid droplet on a cylinder surface. The main fluid point cloud is shown with spherical glyphs: blue for interior points, and red for free boundary points. The background point cloud for the shape of the overall computational domain is shown with black points. To enhance visualization, we have added a grey shell by tessellating the background point cloud.}
	\label{Fig:CylinderFreeMove_IC}%
\end{figure}

~\\
Consider a cylinder with axis along the $x$ direction: $y^2 + z^2 = r^2$, limited by $-0.5 \leq x \leq 0.5$. Here, we set $r=0.5$. An initial fluid droplet is taken as the intersection of the cylinder with the sphere $x^2 + y^2 + (z-0.5)^2 \leq 1$, as shown in Figure~\ref{Fig:CylinderFreeMove_IC}.

\begin{figure}
	\centering
	\includegraphics[width=0.5\textwidth]{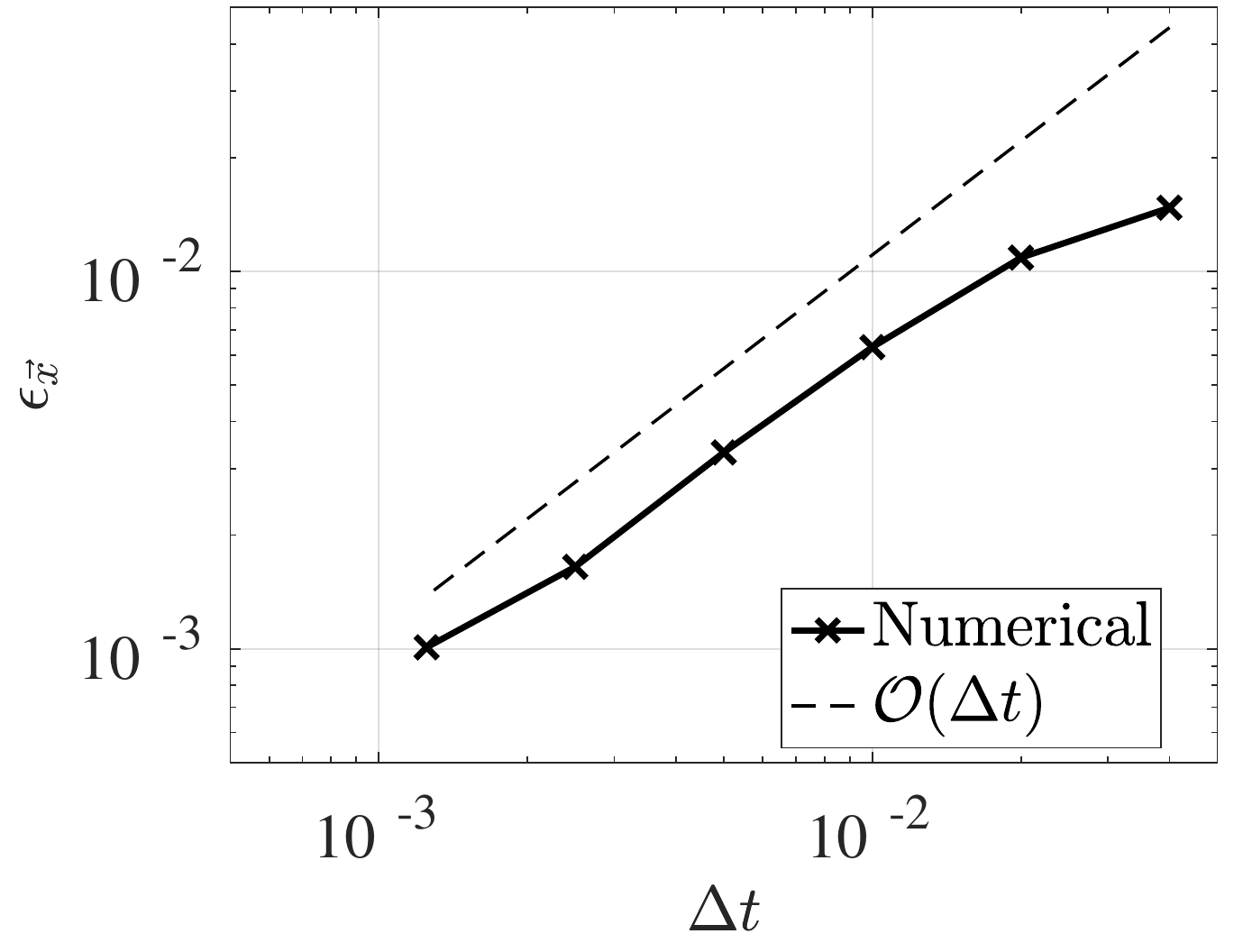} 
	\caption{Error in location of the fluid droplet moving on a fixed cylinder.}
	\label{Fig:CylinderFreeMove_Error}%
\end{figure}

The surface is stationary $\vec{w}=0$, and the droplet moves with a prescribed velocity of 
\begin{equation}
	\vec{v} = ( 0.2\sin( 2\pi t), -\pi z , \pi y )^T\,.
\end{equation}
Note that this velocity is tangential to the surface. The droplet moves in a sinusoidal motion while rotating along the cylinder. At every integer multiple of $t=2$, the droplet returns to its original location. We measure the error of the location of the fluid domain as 
\begin{equation}
	\epsilon_{\vec{x}} = \| \vec{x}_{num} - \vec{x}_{exact} \| \,,
\end{equation}
where $\vec{x}_{num}$ is the location of the geometric center of the fluid point cloud, and $\vec{x}_{exact}$ is the analytical path travelled by the center of the initial fluid domain. For coarse point clouds with $N=138$, and $N^{\text{back}}=1574$ points initially, the convergence of errors at $t=2$ with a decreasing time step is shown in Figure~\ref{Fig:CylinderFreeMove_Error}. 

Only linear convergence with $\Delta t$ is seen. This could have multiple reasons. Firstly, the splitting of the movement is of first oder. Secondly, as also observed in \cite{Suchde2019_MovingSurfaces}, the largest error in the movement comes at the very first time step, where a first order move has to be done since no previous velocity is available. 

\subsubsection{Moving droplet on a moving cylinder}
\begin{figure}
  \centering
  	\includegraphics[width=0.05\textwidth]{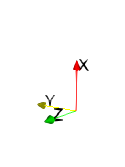} 
	\includegraphics[width=0.22\textwidth]{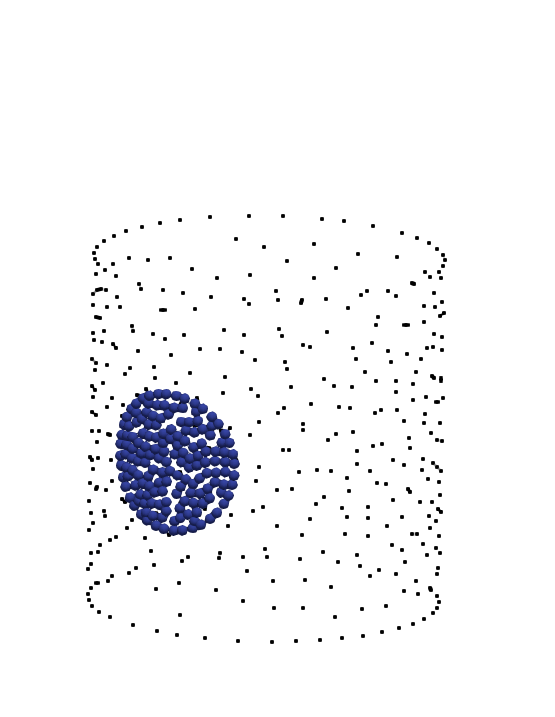} 
	\includegraphics[width=0.22\textwidth]{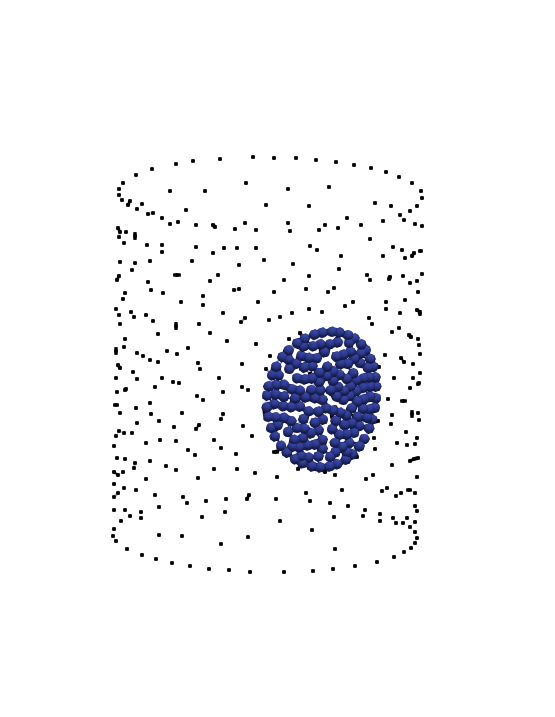} 
	\includegraphics[width=0.22\textwidth]{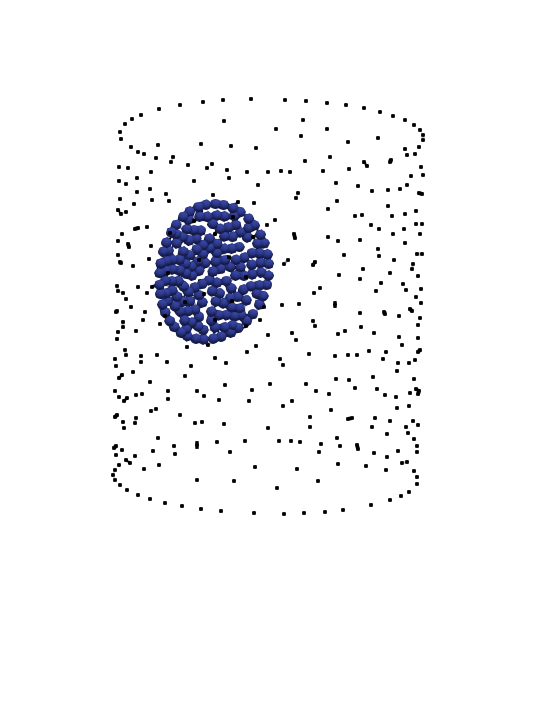} 	
	\includegraphics[width=0.22\textwidth]{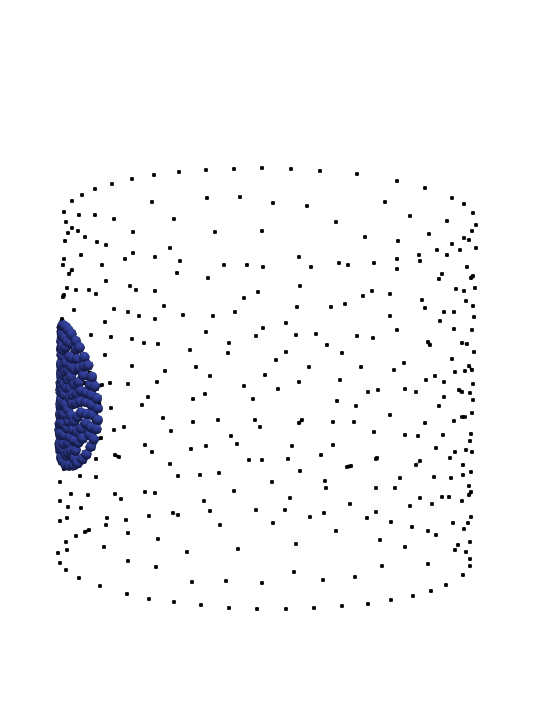} 
  \caption{Evolution of the surface domain and the fluid for the moving droplet on a moving cylinder test case. $t=0$, $t=0.32$, $t=1.32$, and $t=1.76$. The fluid is shown with blue spherical glyphs. And the surface is shown with black points. Note that the scale and coordinate limits are the same in each of the images.}
  \label{Fig:MovingDropMovingCylinder}%
\end{figure}

~\\
We now extend the above example to a moving cylinder. The same initial cylinder and fluid droplet is considered as above. In this case, the cylinder is moving and deforming with
\begin{equation}
	\vec{w} = \left( 0.5  \sin(\pi  t) , 0  ,0  \right)^T 
			- 0.5 \cos(2 \pi  t) \vec{n}  \,,
\end{equation}
where $\vec{n}$ is the outward pointing unit normal of the cylinder surface. The first term causes a translation motion of the cylinder in the direction of the axis, while the second causes the cylinder to expand and contract. The normal component of the surface velocity can be integrated to obtain the analytical radius of the cylinder $r(t) = 0.5 - 0.5\frac{\sin(2\pi t)}{2\pi}$. The same velocity $\vec{v}$ is taken as in the previous example. Note that this velocity still remains tangential to the surface throughout the temporal domain. Snapshots of the evolution of the domain are shown in Figure~\ref{Fig:MovingDropMovingCylinder}.

\begin{figure}
	\centering
	\includegraphics[width=0.5\textwidth]{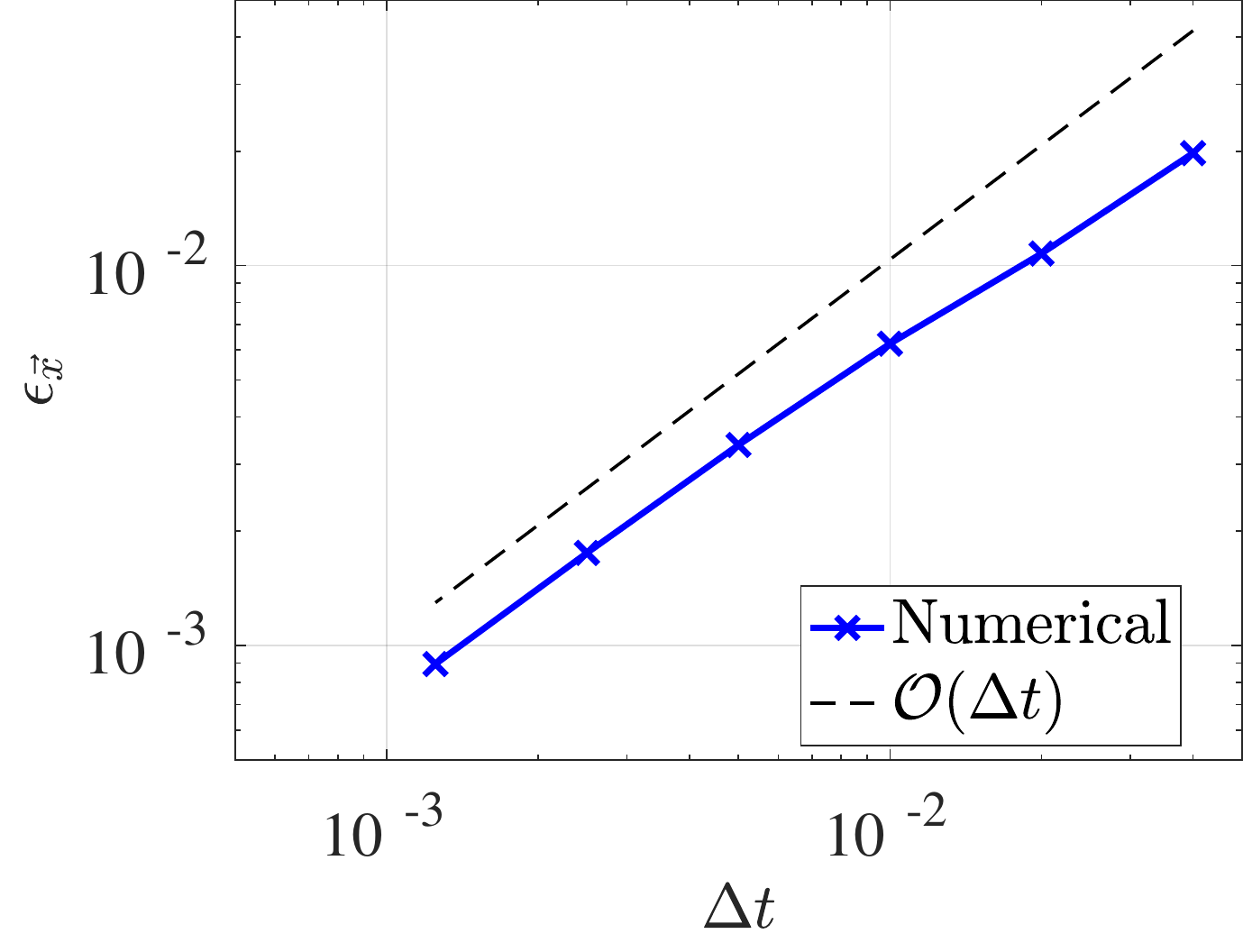} 
	\caption{Error in location of the fluid droplet moving on a moving cylinder. Errors observed are very similar to those of the fixed cylinder example, shown in Figure~\ref{Fig:CylinderFreeMove_Error}.}
	\label{Fig:CylinderFreeMove_Error2}%
\end{figure}

The same initial point clouds are used as in the previous example. At $t=2$, the cylinder domain, and the fluid droplet both return to its initial shape and position. Errors are thus measured as in the previous example, and are shown in Figure~\ref{Fig:CylinderFreeMove_Error2}. The errors do not differ greatly from the fixed cylinder case of the previous example.

\subsubsection{Movement with distortion}~\\
We now consider a case where the fluid domain not just deforms, but also undergoes a topological change. The deformation of the surface itself illustrates the benefit of a background point cloud, rather than a background mesh.

\begin{figure}
  \centering
  \includegraphics[width=0.3\textwidth]{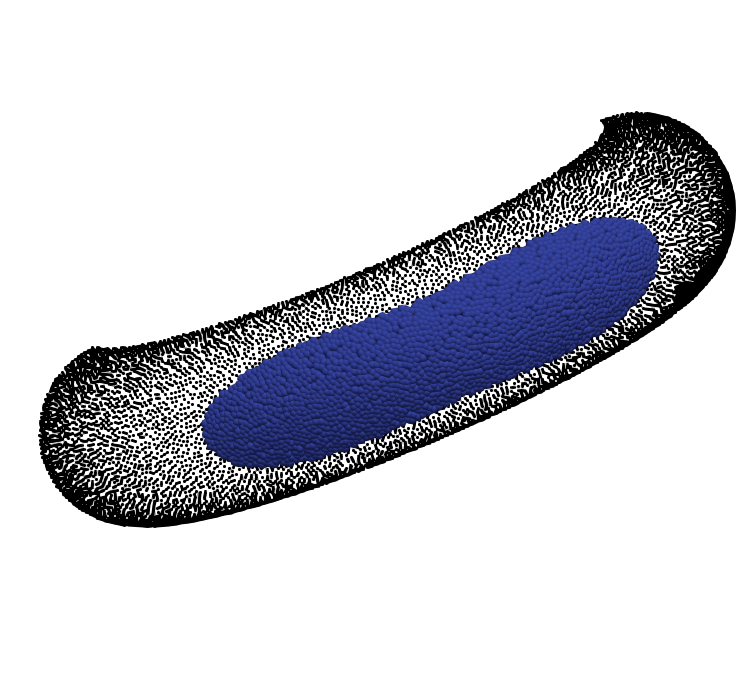}
  \includegraphics[width=0.3\textwidth]{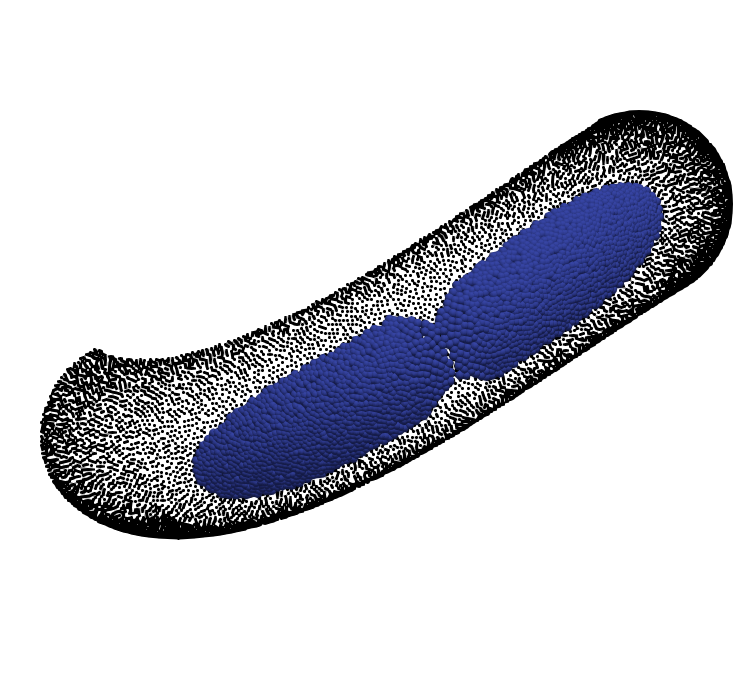}
  \includegraphics[width=0.3\textwidth]{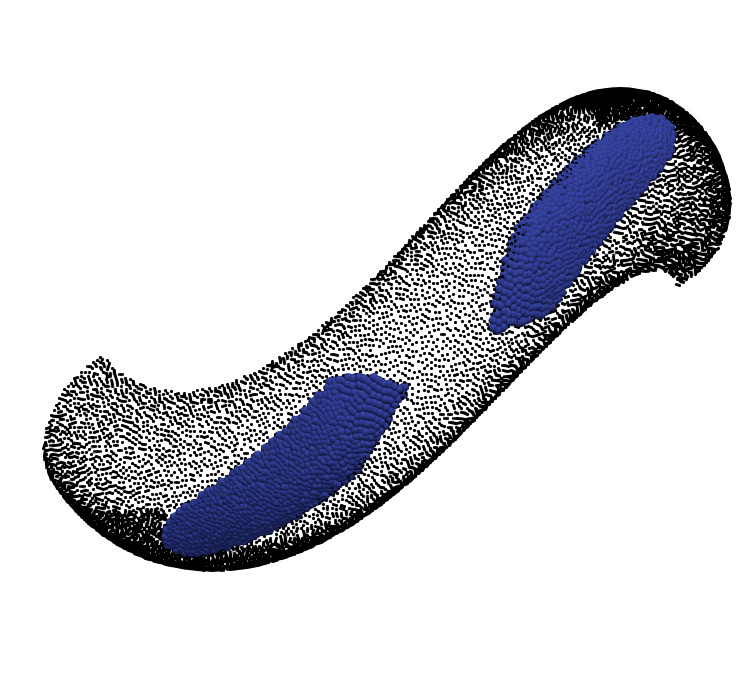}\\
  \includegraphics[width=0.05\textwidth]{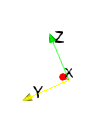}\\
  \includegraphics[width=0.3\textwidth]{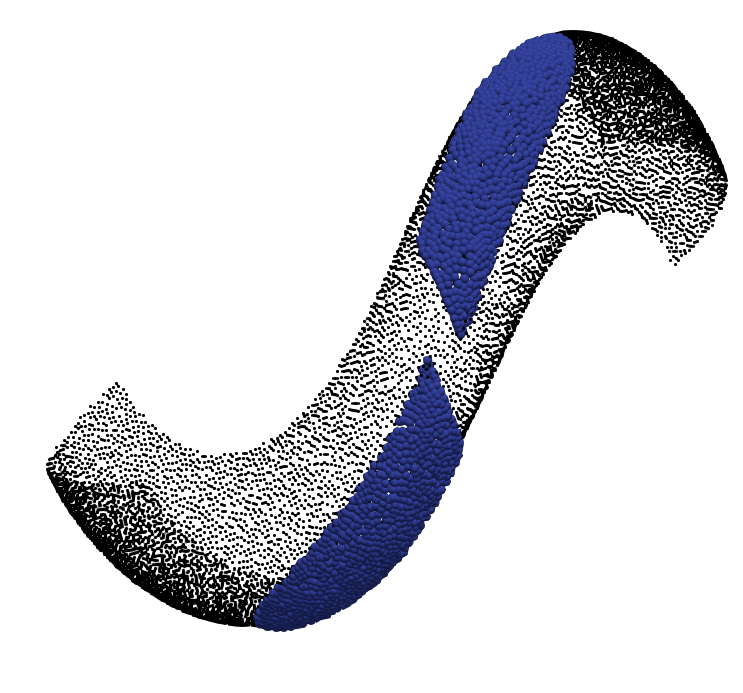}
  \includegraphics[width=0.3\textwidth]{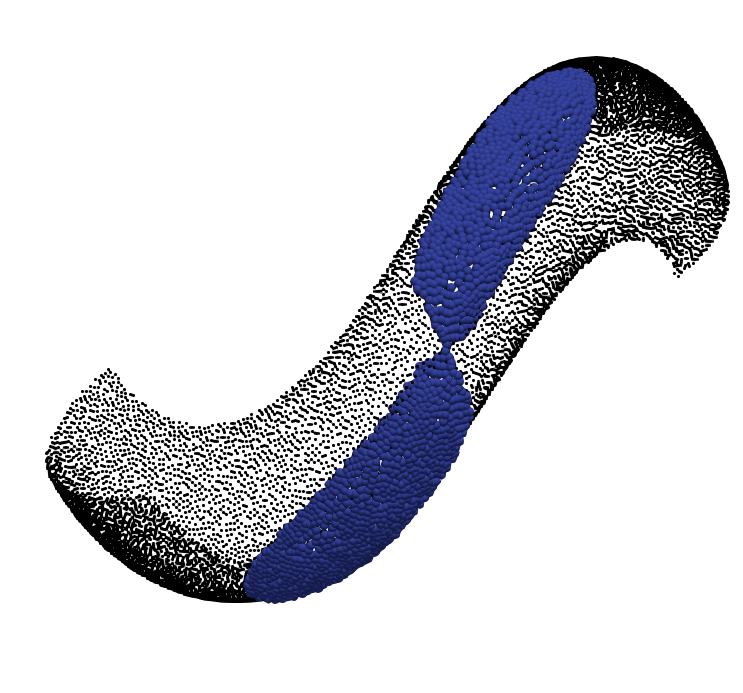}  
  \includegraphics[width=0.3\textwidth]{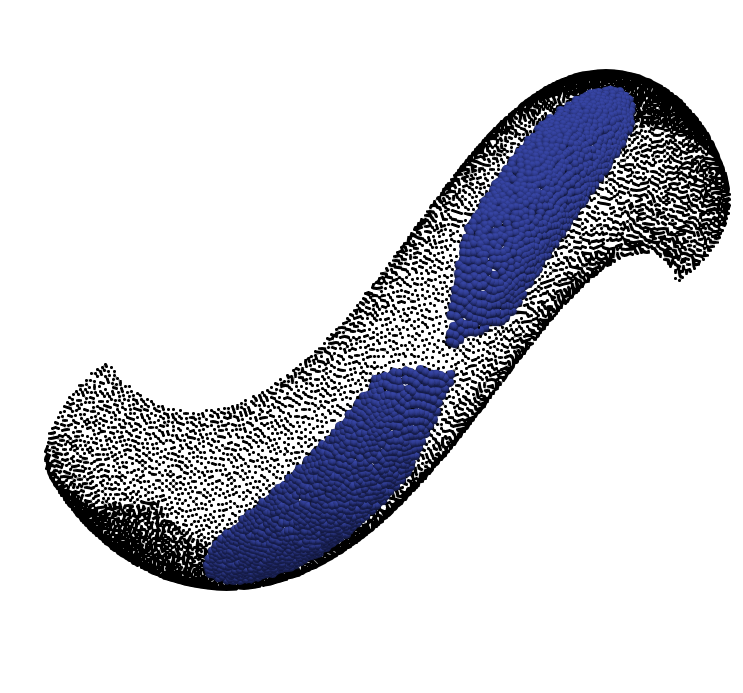}
  \caption{Fluid droplet undergoing a change of topology on a distorting surface. Clockwise from top left: $t=0$~(top left), $t=0.18$~(top center), $t=0.525$~(top right), $t=0.72$~(bottom right), $t=0.855$~(bottom center), and  $t=1.14$~(bottom left). The fluid is shown with blue spherical glyphs, and the surface is shown with black points.}
  \label{Fig:DropTorusTwist_Video}%
\end{figure}

For the surface movement, we consider a case similar to that proposed in \cite{Suchde2019_MovingSurfaces}. The surface domain is a quarter of a torus with major radius $0.75$ and minor radius $0.25$, as shown in the first image of Figure~\ref{Fig:DropTorusTwist_Video}. The surface velocity is given by
\begin{equation}
	\vec{w} = \left( -yz, 0, xz \right)^T \,,
\end{equation}
which causes an overall twisting motion in two different directions. A large fluid droplet is considered at the center of the surface, as shown in the first image of Figure~\ref{Fig:DropTorusTwist_Video}. The fluid velocity is given by
\begin{equation}
	\vec{v} = \boldsymbol{P} (0,\text{sign}(y) \sin(2\pi t),0) ^T \,,
\end{equation}
where $\boldsymbol{P}$ projects the vectors that it acts upon to the tangent space of the surface domain. A more detailed explanation of $\boldsymbol{P}$ is given in the next section. This velocity causes the droplet to split into two parts as it starts moving. 

Finer point clouds are used in comparison to the previous examples to accurately capture the deformation. The simulations start with $N=1902$ points and $N^{\text{back}}=10\,499$ points. The evolution of the fluid and surface domains are shown in Figure~\ref{Fig:DropTorusTwist_Video}. This example illustrates the capability of the present method to handle complex geometries and arbitrary motion.
%


\section{Surface Derivatives}
\label{sec:DiffOperators}

In this section, we present a concise overview of operators on surfaces \cite{Fries2018, Suchde2019_StaticSurfaces} and numerical approximations of surface differential operators using a meshfree Generalized Finite Difference Method~(GFDM) \cite{Suchde2019_StaticSurfaces}. 

\subsection{Differential Operators on Surfaces}

Consider an open subset $\Omega$ of $\mathbb{R}^3$ containing the surface $M$. Thus, we have $M \subset \Omega \subset \mathbb{R}^3$. For a scalar-valued function defined on the surface $f: M \rightarrow \mathbb{R}$, an extension of $f$ to $\Omega$ is defined as $\hat{f}: \Omega \rightarrow \mathbb{R}$ such that $f$ and $\hat{f}$ agree on $M$. Note that the extension is not unique. The so-called tangential differential calculus framework provides a method to define the derivatives of $f$ along the surface $M$ in relation to the conventional derivatives of $\hat{f}$ on $\Omega$. This definition is independent of which extension is used, with the assumption that only ``smooth enough" extensions are considered. In this framework, the surface gradient of $f$ is defined as
\begin{equation}
	\label{Eq:SurfGrad}
	\nabla_M f = \boldsymbol{P} \nabla \hat{f}\,,
\end{equation}
with the projection operator (to the tangent plane)
\begin{equation}
	\boldsymbol{P} = \boldsymbol{I} - \vec{n}\vec{n}^T\,.
\end{equation}
For a vector valued function $\vec{v}: M \rightarrow \mathbb{R}^3$, its extension $\vec{\hat{v}}: \Omega \rightarrow \mathbb{R}^3$ is defined component wise. Thus, if $\vec{v}=(u,v,w)^T$, then $\vec{\hat{v}} = (\hat{u}, \hat{v}, \hat{w} )^T$. The surface divergence of $\vec{v}$ is defined as
\begin{equation}
	\label{Eq:SurfDiv}
	\nabla_{M}\cdot \vec{v} =  \left(\boldsymbol{P}\nabla\right) \cdot \vec{\hat{v}}\;.
\end{equation}
Combining Eqs.\,\eqref{Eq:SurfGrad} and \eqref{Eq:SurfDiv} gives the definition of the surface Laplacian or Laplace-Beltrami of $f$
\begin{align}
	\Delta_{M} f &= \nabla_{M}\cdot \nabla_{M} f\,,\\
	&= \left(\boldsymbol{P}\nabla \right) \cdot \left( \boldsymbol{P} \nabla \hat{f} \right) \label{Eq:SurfLap}\,.
\end{align}
The directional surface gradient of a vector valued function is given by
\begin{align}
	\nabla_M^{\text{dir}} \vec{v} &= \left( \nabla \vec{\hat{v}} \right) \boldsymbol{P} \,, \\ 
	&= \left( \nabla_M u \;, \nabla_M v \;, \nabla_M w \right) ^T \,.
\end{align}
Furthermore, the covariant surface gradient of a vector valued function is defined as the projection of the directional gradient to the tangent plane, and is given by
\begin{align}
	\label{Eq:CovVecGradient}
	\nabla_M^{\text{cov}} \vec{v} &= \boldsymbol{P}  \nabla_M^{\text{dir}} \vec{v}\,,\\
	&=  \boldsymbol{P}  \nabla \vec{\hat{v}} \boldsymbol{P} \,.
\end{align}
For the divergence of tensor-valued functions, we follow the convention of column-wise divergences. 

The advantages of this framework of defining surface operators are that it avoids differential geometry complexities and avoids the need to parametrize the surface. Furthermore, it provides a straight forward interpretation for strong form PDE solvers, which shall be used in the present work. This comes at the cost of having to solve an extra unknown for vector fields. For example, if the surface can be parametrized, the velocity field only needs two components to describe it entirely. However, solving it in a global Eucledian setting requires three velocity components. For more details of the tangential differential calculus framework, we refer to \cite{Fries2018, Suchde2019_StaticSurfaces}. 

\subsection{Numerical Surface Derivatives}
\label{sec:NumDO}

Generalized Finite Difference Methods~(GFDM) are strong form collocation methods to discretize differential operators. They have been shown to be robust methods, and have been widely applied to solve PDEs on volume domains (for example, \cite{Drumm2008, Fan2018,Gavete2017,Katz2010,Luo2016} ). In the present context, we use GFDMs to discretize the surface differential operators. The general idea of computing numerical derivatives that we use here is to project neighbouring points to a tangent space, compute regular volume 2D derivatives in the tangent space, and rotate them back to the surface to obtain surface derivative stencils. This notion was introduced in \cite{Demanet2006} for the surface Laplacian on meshes, and was extended in our earlier work \cite{Suchde2019_StaticSurfaces} for the use on general differential operators on point clouds. For the sake of completeness, we give an brief overview below.

The normal extension of a function is the extension which satisfies $\vec{n}\cdot\nabla \hat{f} \equiv 0$. For normally extended functions, Eqs.\,\eqref{Eq:SurfGrad}, \eqref{Eq:SurfDiv} and \eqref{Eq:SurfLap} reduce to (see \cite{Suchde2019_StaticSurfaces})
\begin{align}
	\nabla_M f &= \nabla \hat{f}\,,\\
	\nabla_{M}\cdot \vec{v} &=  \nabla \cdot \vec{\hat{v}}\,,\\
	\Delta_{M} f &= \Delta \hat{f}\,.
\end{align}
Thus, the surface gradient, divergence and Laplacian can be computed by approximating the corresponding conventional derivative of the normal extension. We use this to numerically compute the surface derivatives. 

Local normal extensions are used to project neighbouring points to the tangent space. For the surface gradient, at a numerical point $i$ at location $\vec{x}_i$ on the surface, we compute the tangential components of $\nabla \hat{f}$ on the tangential plane $T_i$ spanned by the tangents $\vec{t}_{1,i}$ and $\vec{t}_{2,i}$. For this, the neighbouring points $j \in S_i$ are projected along the normal $\vec{n}_i$ to the locations $j_{T_i}$. First, we approximate the numerical gradient on the tangent plane $\nabla_T$
\begin{equation}
	\nabla_{T} \hat{f} \approx \widetilde{\nabla}_{T} \hat{f} = 
			\left( 	\sum_{j\in S_i} c_{ij_T}^{t_1}\hat{f}_{j_T}  \;\; , \;\;
					\sum_{j\in S_i}	c_{ij_T}^{t_2}\hat{f}_{j_T}  \right)^T \,,
\end{equation}
where the overhead $\sim$ indicates the numerical operator. Stencil coefficients $c_{ij_T}$ are computed using a GFDM approach \cite{Gavete2017}, which employs a least squares procedure shown below. Monomials up to a certain order are differentiated exactly
\begin{align}
	\sum_{j\in S_i}c_{ij_{T}}^{t_k}m_{j_{T}} &= \frac{\partial}{\partial t_k} m (\vec{x}_i) \qquad \forall m\in\mathcal{P}_{T}\,,\label{Eq:TP_Consistency}\\
	\text{min } J_i &= \sum_{j\in S_i} \left( \frac{c_{ij_{T}}^{t_k}}{W_{ij_{T}}} \right)^2\,, \label{Eq:TP_Min}
\end{align}
for $k=1,2$, where $\mathcal{P}_{T}$ is the set of monomials, taken up to order $2$ in the present work, in $\vec{t}_{1,i}$ and $\vec{t}_{2,i}$ on the tangent plane. Now, the numerical surface gradient is defined as
\begin{equation}
	\label{Eq:SurfGrad_FullDefinition}
	\widetilde\nabla_{M,i} f = %
					\left( 
							\sum_{j \in S_i} c_{ij}^{M,x} f_j \;\; , \;\;
							\sum_{j \in S_i} c_{ij}^{M,y} f_j \;\; , \;\;	
							\sum_{j \in S_i} c_{ij}^{M,z} f_j \right)^T \,,		
\end{equation}
where $c_{ij}^{M,x}$ are the stencil coefficients for the surface gradient in the $x$ direction, and similarly for the other directions. These stencil coefficients are given by 
\begin{equation}
	\left(  	c_{ij}^{M,x} \;\; , \;\;
			c_{ij}^{M,y} \;\; , \;\;
			c_{ij}^{M,z} 
	\right)^T = %
	R^T 		\left( 
			c_{ij}^{t_1} \;\; , \;\;												   									c_{ij}^{t_2} \;\; , \;\;
			c_{ij}^{n} 
		\right)^T\,,
\end{equation}
with $c_{ij}^n = 0$ and a rotation matrix $R$
\begin{equation}
	R^T = \left(\begin{array}{ccc} 
	\vec{t}_1 & \vec{t}_2 & \vec{n} 
	\end{array}\right)\,,
\end{equation}
The numerical gradient stencils also provide the numerical divergence. A similar procedure is followed for the numerical surface Laplacian defined as
\begin{equation}
	\label{Eq:LapBeltNumericalDefinition}
	\widetilde{\Delta}_{M} f = \sum_{j\in S_i} c_{ij}^{\Delta_M}f_j\,.
\end{equation}
Rotational invariance of the Laplace operator means that the surface Laplacian is directly given by the tangent place Laplacian $c_{ij}^{\Delta_M} = c_{ij}^{\Delta_T}\;$, which is computed by
\begin{align}
	\sum_{j\in S_i}c_{ij_{T}}^{\Delta_{T}}m_{j_{T}} &= \Delta_{T} m (\vec{x}_i) \qquad \forall m\in\mathcal{P}_{T}\,,\label{Eq:TPL_Consistency}\\
	\text{min } J_i &= \sum_{j\in S_i} \left( \frac{ c_{ij_{T}}^{\Delta_{T}} } {W_{ij_{T}} } \right)^2\,. \label{Eq:TPL_Min}
\end{align}
We do not discretize the covariant vector derivative directly. Only surface gradients and Laplacians are the computed numerical differential operators. For vector fields, the application of the numerical gradient operator gives us the directional vector gradient. The numerical covariant vector gradient can then be computed using the numerical projection operator, according to Eq.\,\eqref{Eq:CovVecGradient}. For more details of both the theory and implementation, including the necessary stabilization of the surface Laplacian defined above, we refer to \cite{Suchde2019_StaticSurfaces}.


\section{Navier--Stokes Equations on Manifolds}
\label{sec:NSman}

We consider the incompressible Navier--Stokes equations posed on a curved surface \cite{Fries2018, Jankuhn2017, Miura2017}. The surface stress tensor is given by
\begin{equation}
	\boldsymbol{S} = 	\eta \left[ \left(\nabla_M^{\text{cov}} \vec{v}\right) + \left( \nabla_M^{\text{cov}} \vec{v} \right)^T \right] \,,
\end{equation}
for (dynamic) viscosity $\eta$. Using this, the incompressible surface Navier--Stokes equations read as
\begin{align}
	\frac{\partial\vec{v}}{\partial t} + \left( \vec{v}\cdot \nabla_M^{\text{cov}} \right) \vec{v} + \left( \vec{w}\cdot\nabla \right) \vec{v} &= \frac{1}{\rho} \boldsymbol{P} \nabla_M \cdot  \boldsymbol{S} - \frac{1}{\rho}\nabla_M p + \vec{g} \label{Eq:CoMo}\,,\\
	\nabla_M \cdot \vec{v} &= 0 \,,\label{Eq:CoM}\\
	\vec{v}\cdot \vec{n} &= 0 \label{Eq:vTan}\,,
\end{align}
with fluid velocity $\vec{v}$ on the surface, surface velocity $\vec{w}$, pressure $p$, and gravity and body forces term $\vec{g}$. Eq.\,\eqref{Eq:CoMo} arises from momentum conservation, Eq.\,\eqref{Eq:CoM} is the incompressibility constraint, and Eq.\,\eqref{Eq:vTan} indicates that the velocity field lies entirely in the tangent bundle, with no normal component. To simplify the model, inertial source terms arising from the movement of the surface are ignored here. 


\subsection{Numerical scheme}
\label{sec:NumericalScheme}

We follow a Chorin projection method \cite{Chorin1968} and modify it for the needs of the surface Navier--Stokes equations. First, we split the viscous term into the projection of the component wise Laplacian, which resembles the conventional viscous term in volumetric Navier--Stokes equations, and extra rotational components
\begin{equation}
	\label{Eq:DSTsplit} 
	\frac{1}{\rho} \boldsymbol{P} \nabla_M \cdot  \boldsymbol{S} \left( \vec{v} \right) = 
		\left( %
		\frac{\eta}{\rho}\boldsymbol{P}\Delta_M \vec{v} \right) +
		\left( %
		\frac{1}{\rho} \boldsymbol{P} \nabla_M \cdot  \boldsymbol{S} ( \vec{v} ) - 
		\frac{\eta}{\rho}\boldsymbol{P}\Delta_M \vec{v} \right) \,.
\end{equation}
Note that we assume $\eta$ to be constant in space. We start the scheme with a Lagrangian movement step, Eq.\,\eqref{Eq:MoveOverall}, as explained in Sections~\ref{sec:LagMove} and \ref{sec:FreeBoundary}. This is followed by semi-implicit intermediate velocity solve, which treats the first term in Eq.\,\eqref{Eq:DSTsplit} implicitly and second one explicitly
\begin{equation}
\label{Eq:Vstar}
	\frac{\vec{v}^{\,*} - \vec{v}^{\,(n)} }{\Delta t} - 
	\frac{\eta}{\rho} \boldsymbol{P} \Delta_M \vec{v}^{\,*} =
	\frac{1}{\rho} \boldsymbol{P} \nabla_M \cdot  \boldsymbol{S} ( \vec{v}^{\,(n)} ) - 
	\frac{\eta}{\rho} \boldsymbol{P} \Delta_M \vec{v}^{\,(n)}  - 
	\frac{1}{\rho}\nabla_M p^{(n)} + 
	\vec{g} \,.
\end{equation}	
In projection methods for volume domain flows, this intermediate velocity would then be projected to a divergence free space with the help of a correction pressure. Here, we first project the intermediate velocity to the tangent plane
\begin{equation}
\label{Eq:Vstarstar}
	\vec{v}^{\,**} = \boldsymbol{P}  \vec{v}^{\,*} \,.
\end{equation}
This velocity is now projected to a surface divergence free space
\begin{equation}
\label{Eq:Vupdate}
	\vec{v}^{\,(n+1)} = \vec{v}^{\,**} - \frac{\Delta t}{\rho}\nabla_M p_{\text{corr}} \,,
\end{equation}	
where the correction pressure $p_{\text{corr}}$ is first computed by applying the surface divergence operator to Eq.\,\eqref{Eq:Vupdate}, to give the pressure Poisson equation
\begin{equation}
\label{Eq:Ppoisson}
	\frac{\Delta t}{\rho} \Delta_M p_{\text{corr}} = \nabla_M \cdot \vec{v}^{\,**} \,,
\end{equation}	
for spatially constant $\rho$. Finally, the pressure is updated 
\begin{equation}
\label{Eq:Pupdate}
	p^{(n+1)} =  p^{(n)} + p_{\text{corr}} \,.
\end{equation}
%

We note that the above scheme automatically enforces Eq.\,\eqref{Eq:vTan}. Firstly, $\vec{v}^{\,**} \cdot \vec{n} = 0$ by Eq.\,\eqref{Eq:Vstarstar}. Furthermore, we know that the surface gradient of a scalar-valued function always lies in the tangential plane, i.e., $\vec{n} \cdot \nabla_M \phi = 0$, $\forall \phi: M \rightarrow \mathbb{R}$. Thus, $\vec{v}^{\,(n+1)}$ according to Eq.\,\eqref{Eq:Vupdate} also lies in the tangential plane as it is given by a linear combination of $\vec{v}^{\,**}$ and $\nabla_M p_{\text{corr}}$.

To obtain the discrete systems, the continuous operators are replaced with the corresponding discrete ones as defined in Section~\ref{sec:NumDO}, with the appropriate boundary conditions. So, for example, the discrete version of the pressure Poisson system Eq.\,\eqref{Eq:Ppoisson} on a closed surface (without boundaries) reads as
\begin{equation}
\label{Eq:DiscPpoisson}
	\frac{\Delta t}{\rho} \sum_{j \in S_i} c_{ij}^{\Delta_{M}} p_{\text{corr}, j} = \sum_{j \in S_i} c_{ij}^{M,k} v^{\,**}_k \,,
\end{equation}
for $k = x,y,z$, and the velocity $\vec{v}^{\,**} = (v^{\,**}_x, v^{\,**}_y, v^{\,**}_z)^T$, and $i=1,2,\dots,N$.




\subsection{Boundaries and Boundary Conditions}
\label{sec:BC}

One of the advantages of meshfree GFDMs over other meshfree methods such as SPH is the ease of enforcing a variety of boundary conditions. This advantage also carries over to the surface PDE context, as illustrated in \cite{Suchde2019_StaticSurfaces}. For example, to enforce boundary conditions on the pressure Poisson equation, the corresponding row of Eq.\,\eqref{Eq:DiscPpoisson} is simply replaced with the appropriate discrete boundary condition. 

For the Lagrangian framework introduced earlier, a few extra modifications are needed for the imposition of certain boundary conditions. Firstly, points on inflow and outflow boundaries are kept fixed with respect to the fluid velocity $\vec{v}$, and are only moved with the surface velocity $\vec{w}$. This ensures a certain regularity of the point cloud at the boundary. Moreover, if points at the inflow were moved in a Lagrangian sense, new points would have to be added at the inflow boundary. This would require an extrapolation of field values there, since the location would be outside the moved point cloud. On the other hand, if the inflow boundary points are fixed, points can be added in the resultant hole formed in the inflow direction (see Section~\ref{sec:Distortion}), and those points would require an interpolation of field values with points present on all sides, and not an extrapolation. Furthermore, points crossing an outflow boundary are deleted. Unlike inflow and outflow boundaries, points on slip boundaries undergo the Lagrangian motion even with $\vec{v}$, but with the extra constraint that they must move only along the boundary. During the merging procedure explained in Section~\ref{sec:Distortion}, priority is given to all boundary points, which ensures that an interior point and a boundary point are not merged, which would have led to artificially distorting the boundary.

Neumann boundary conditions for the pressure field $\vec{\nu}\cdot \nabla p = r$ are discretized as
\begin{equation}
	\sum_{j \in S_i} c_{ij}^{M,\nu} p_j = r \,,
\end{equation}
with
\begin{equation}
	c_{ij}^{M,\nu} = \vec{\nu} \cdot \left( c_{ij}^{M,x}, c_{ij}^{M,y}, c_{ij}^{M,z} \right) ^T \,.
\end{equation}

For free boundaries, the boundary conditions are based on surface stresses
\begin{align}
		\vec{t}^T \, \boldsymbol{S}(\vec{v}) \, \vec{\nu} &= F_{ext} \,, \label{Eq:FB_bc1}\\
		\vec{\nu}^T \, \boldsymbol{S}(\vec{v}) \, \vec{\nu} &= p - F_{\sigma} - F_c \,, \label{Eq:FB_bc2}\\
		\vec{n}^T \, \boldsymbol{S}(\vec{v}) \, \vec{\nu} &= 0 \,, \label{Eq:FB_bc3}
\end{align}	
where $F_{ext}$ is the contribution of external forces in the tangential direction, $F_{\sigma}$ represents interface tension effects, and $F_c$ represents contact angle effects. For the sake of simplicity, we only consider the case $F_{ext} = F_{\sigma} = F_c = 0$ in the examples that follow. Since Eq.\,\eqref{Eq:FB_bc3} holds trivially for a tangential velocity field, we replace it with Eq.\,\eqref{Eq:vTan}. 


\section{Validation and Numerical Results}
\label{sec:Results}

In Section~\ref{sec:LagNum}, we validated the Lagrangian framework for flow on (moving) surfaces. We now present a few validation examples for the surface Navier--Stokes equations solver presented above. An irregular point spacing is used in all numerical simulations. Linear systems arising from the semi-implicit intermediate velocity computation Eq.\,\eqref{Eq:Vstar} and the pressure Poisson equation Eq.\,\eqref{Eq:Ppoisson} are both solved using a BiCGSTAB solver \cite{BiCGSTAB}. The number of points used to discretize the domain is always indicated at initial time. Due to the addition and deletion of points to maintain regularity of the point cloud, the number of points will vary slightly in time. The surface Reynolds number is defined in a manner similar to the classical Reynolds number
\begin{equation}
	\text{Re}_{M} = \frac{ \rho U L}{\eta} \,,
\end{equation}
for reference velocity $U$, and length $L$. For the first two examples presented below, which use a manufactured solution, we prescribe a solution for the pressure field $p_{\text{exact}}$, and one for the velocity field $\vec{v}_{\text{exact}}$ that is incompressible, and with $\vec{v}_{\text{exact}} \cdot \vec{n} = 0$. The manufactured gravity term is then computed by substituting these in the momentum equation
\begin{equation}
	\vec{g}_{\text{manufactured}} = 	\frac{\partial\vec{v}_{\text{exact}}}{\partial t} + \left( \nabla_M^{cov} \vec{v}_{\text{exact}}\right) \vec{v}_{\text{exact}} - \frac{1}{\rho} \boldsymbol{P} \nabla_M \cdot  \boldsymbol{S}(\vec{v}_{\text{exact}}) - \frac{1}{\rho}\nabla_M p_{\text{exact}} \,,
\end{equation}
with $\vec{w}=0$ in both cases.


\subsection{Decaying vortices on three quarters of a sphere}

\begin{figure}
  \centering
  \includegraphics[width=0.43\textwidth]{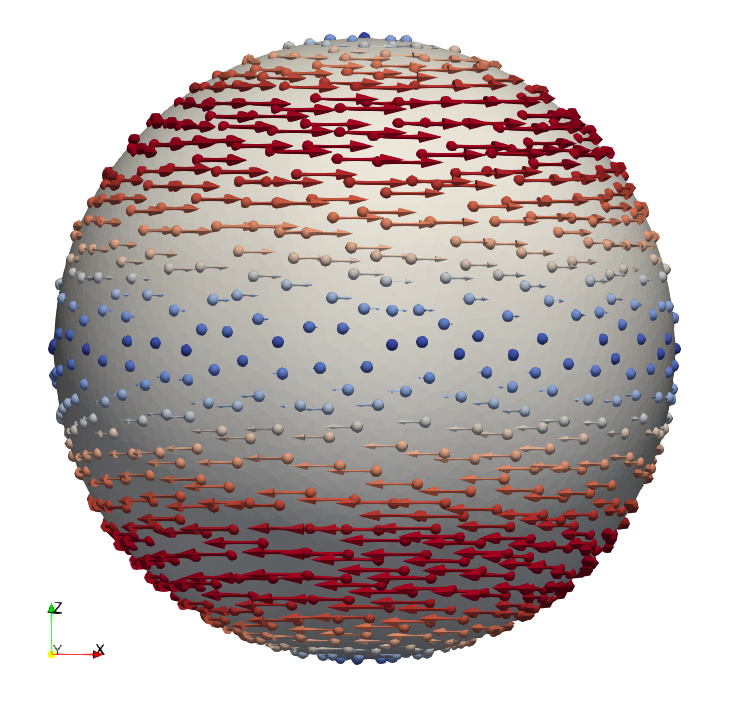}
  \includegraphics[width=0.43\textwidth]{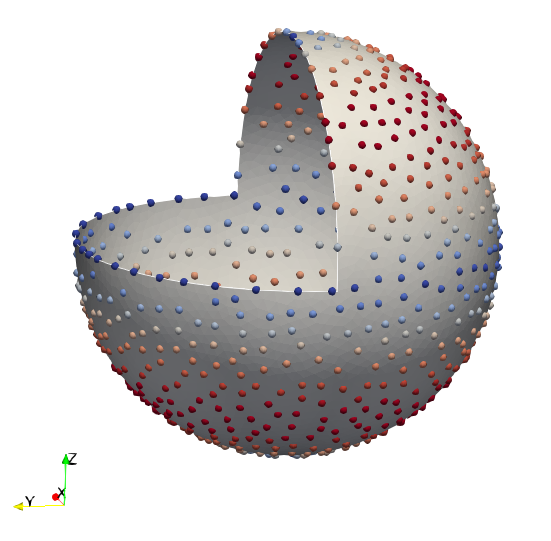}  
  \caption{Decaying vortices on three quarters of a unit sphere. The domain is shown from two different view points. The colour represents the magnitude of the velocity at $t=0$. The arrows in the left figure show the velocity. A grey shell is added to visualize the domain. }
  \label{Fig:Sphere}%
\end{figure}

For the first example, we consider a straight forward case with only Dirichlet boundaries. We start by proposing a manufactured solution on the surface of a unit sphere. The exact solution is given by
\begin{align}
	\vec{v}_{\text{exact}} &= (-yz, xz, 0)^T \exp( -4 \eta t) \,,\\
	p_{\text{exact}} &= -\frac{z^4}{4} \exp(-8 \eta t) \,.
\end{align}
This results in two decaying vortices moving in opposite directions, one each in the positive and negative $z$ hemispheres. The domain is taken to be three quarters of the unit sphere, as shown in Figure~\ref{Fig:Sphere}. Dirichlet boundary conditions are imposed for both the velocity and the pressure on the boundaries, based on the exact solution. For the simulations considered, we set $\rho = 1$ and $\eta = 0.01$. With $U=1$, and $L=1$ is the radius of the sphere, this results in a surface Reynolds number of about $100$.

Errors are measured in a discrete $L^2$ sense, and taken relative to the analytical velocity
\begin{equation}
	\epsilon_2(\vec{v}) = \frac{1}{ \|\vec{v}_{\text{exact}} \|_2}   \left( \frac{1}{N}  \sum_{i=1}^N  \|\vec{v}_i - \vec{v}_{\text{exact}} (\vec{x}_i) \| ^2  \right)^{\frac{1}{2}} \,.
\end{equation}

\begin{figure}
  \centering
  \includegraphics[width=0.43\textwidth]{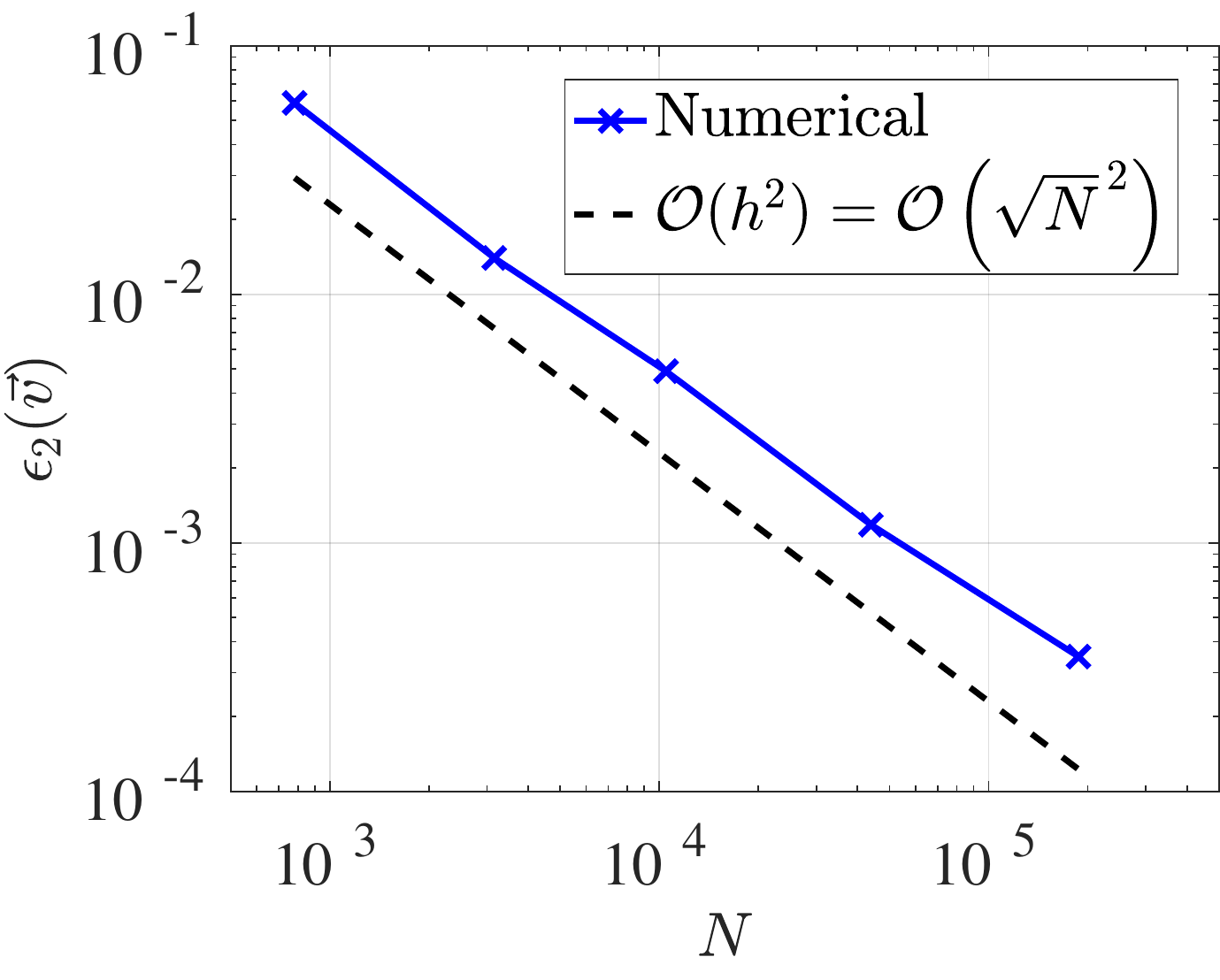}
  \includegraphics[width=0.43\textwidth]{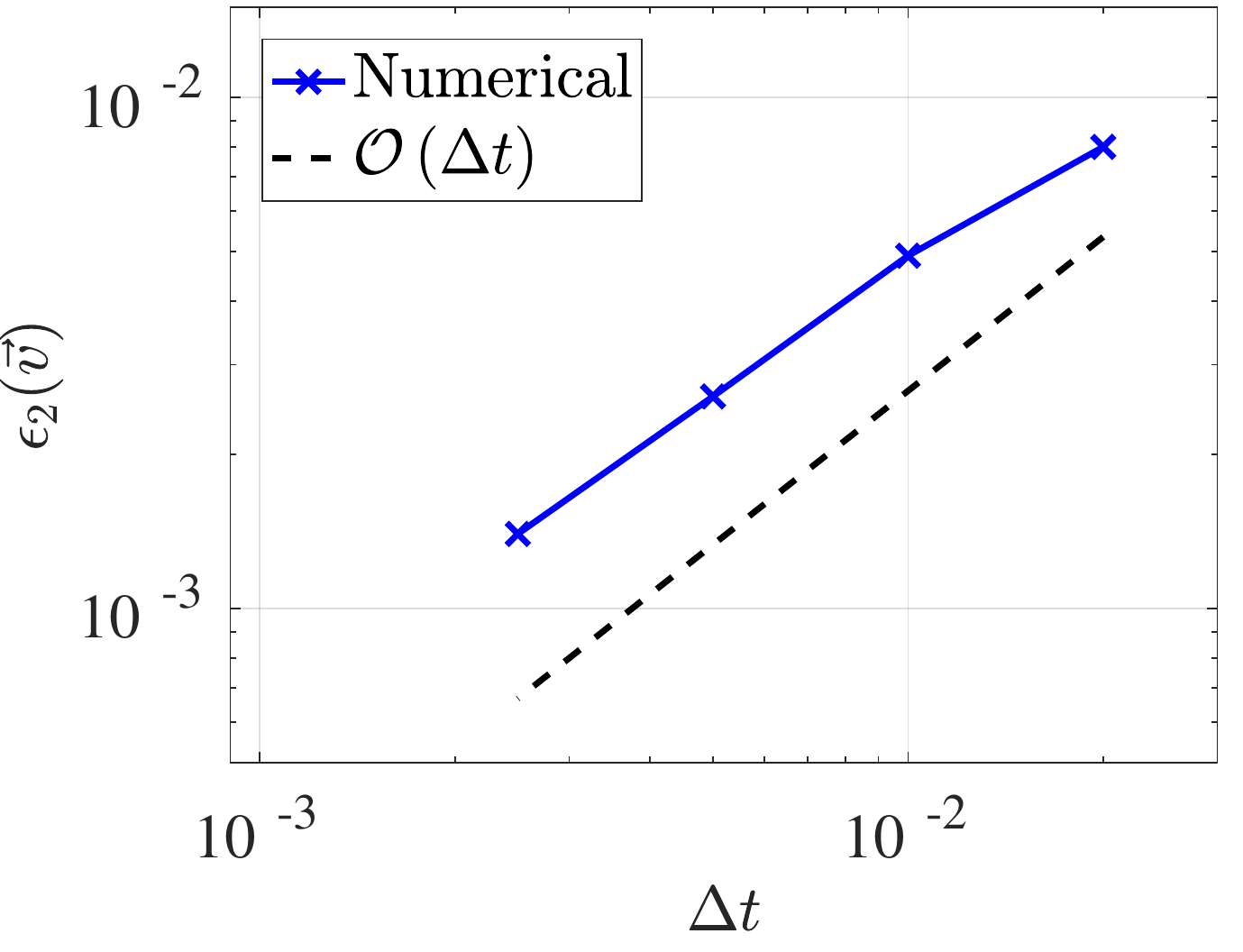}  
  \caption{Decaying vortices on a sphere: Convergence of relative errors with a refining spatial resolution~(left) and a decreasing time step~(right).}
  \label{Fig:Sphere_Convergence}%
\end{figure}

Convergence of errors in the velocity field at $t=1.0$ is shown in Figure~\ref{Fig:Sphere_Convergence} for varying spatial and temporal  resolutions. For the convergence with smoothing length, we use $\Delta t= h^2$,  with the smoothing length consecutively halving from a maximum of $0.4$ to $0.025$. This corresponds to the total number of points in the domain at $t=0$ going from $781$ to $187430$. For the temporal scale convergence, we use $h=0.1$, which corresponds to $N=10488$ at initial time.

Figure~\ref{Fig:Sphere_Convergence} illustrates that an experimental order of convergence of $2$ is observed against varying spatial resolution, which matches the expectation due to the use of up to second order monomials in the numerical differential operators. Furthermore, an approximately first order convergence is observed with varying time step, which matches the use of the first order time integration scheme.


\subsection{Flow on a cylinder}

We now take an example with more involved boundary conditions. Once again, a manufactured solution is used to test the applicability of the method proposed on a larger range of surface Reynolds flow. Consider the surface of a cylinder with axis along the $x$ direction and unit radius, given by $y^2 + z^2 = 1$. The domain is limited by $x \in [0,1]$. The proposed analytical solution is one where both velocity and pressure increase linearly in time
\begin{align}
	\vec{v}_{\text{exact}} &= ( y^2, -z, y)^T t \,,\\
	p_{\text{exact}} &= x t \,.
\end{align}

The boundary at $x=1$ thus acts like an outflow, and that at $x=0$ behaves likes an inflow, and the corresponding treatment for inflow and outflow boundaries respectively are applied, as mentioned in Section~\ref{sec:BC}. Boundary conditions are based on the exact solution. For the pressure, a Dirichlet condition is imposed at the outflow boundary, and a Neumann condition at the inlet. For the velocity field, a Neumann condition is imposed on the outflow outflow with a Dirichlet one at the inlet.

We set $\rho = 1$, and vary $\eta$ to consider different surface Reynolds numbers. For reference values $U =1$, and $L = 1$ for the radius of the cylinder, this gives $\text{Re}_{M} = \frac{1}{\eta}$. Simulations are carried out on a point cloud with $h=0.1$, which corresponds to $N=7119$ points for the initial configuration. For varying $\eta$, the time step considered is given by
\begin{equation}
	\Delta t = \min \{ 0.1, \frac{h^2}{\eta} \} \,.
\end{equation}

Note the small time steps needed for stability for the high viscosity (low surface Reynolds) cases, owing to the explicit terms in the discrete momentum equation. Errors in the velocity field are measured as mentioned in the previous section. A similar definition is used for the pressure field
\begin{equation}
	\epsilon_2(p) = \frac{1}{ \|p_{\text{exact}} \|_2}   \left( \frac{1}{N}  \sum_{i=1}^N  |p_i - p_{\text{exact}} (\vec{x}_i) | ^2  \right)^{\frac{1}{2}} \,.\\
\end{equation}

\begin{figure}
  \centering
  \includegraphics[width=0.45\textwidth]{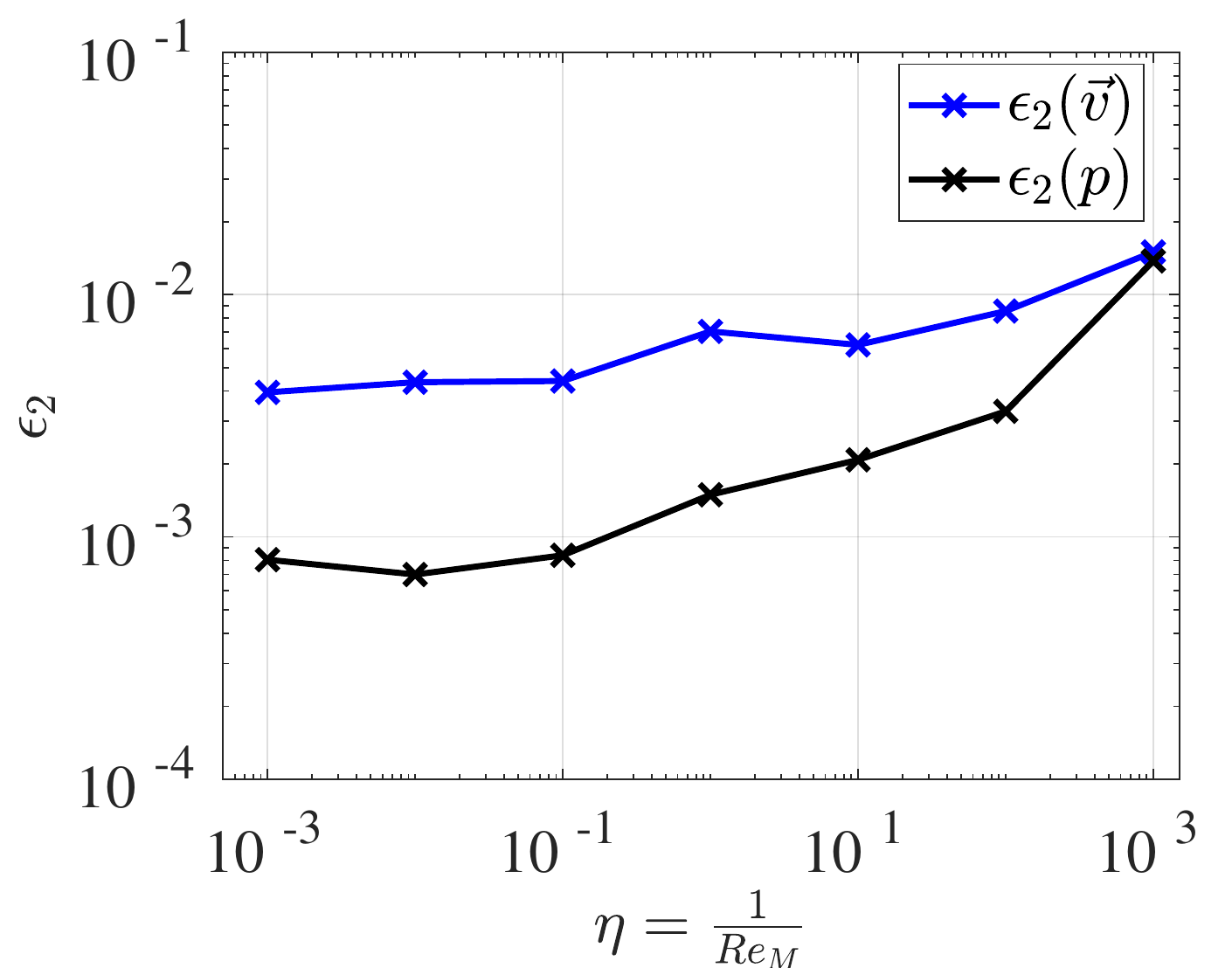}
  \caption{Flow on a cylinder surface: Relative errors in velocity and pressure for different surface Reynolds numbers.}
  \label{Fig:CylinderErrors}%
\end{figure}

The variation of errors in the velocity and pressure field at $t=1.0$ with changing surface Reynolds numbers are shown in Figure~\ref{Fig:CylinderErrors}. A similar error is observed for the velocity field in each of the cases, with a slight increase as the the viscosity increases. However, the increase in error in the pressure field in much more significant. 

It is worth noting here that projection schemes used with strong form meshfree GFDM solvers for volumetric flow often use an additional pressure Poisson equation in each time step to obtain greater accuracy for the pressure field. This is typically done by taking the divergence of the momentum equation with the known velocity field $\vec{v}^{\,(n+1)}$, to obtain the final pressure. Exploring such methods for Navier--Stokes on manifolds, as well as fully implicit velocity schemes and coupled velocity pressure schemes in the meshfree GFDM context remains a topic of investigation beyond the scope of the present paper.

\subsection{Falling Droplet}
\label{sec:DropFall}

\begin{figure}
  \centering
  \includegraphics[width=0.43\textwidth]{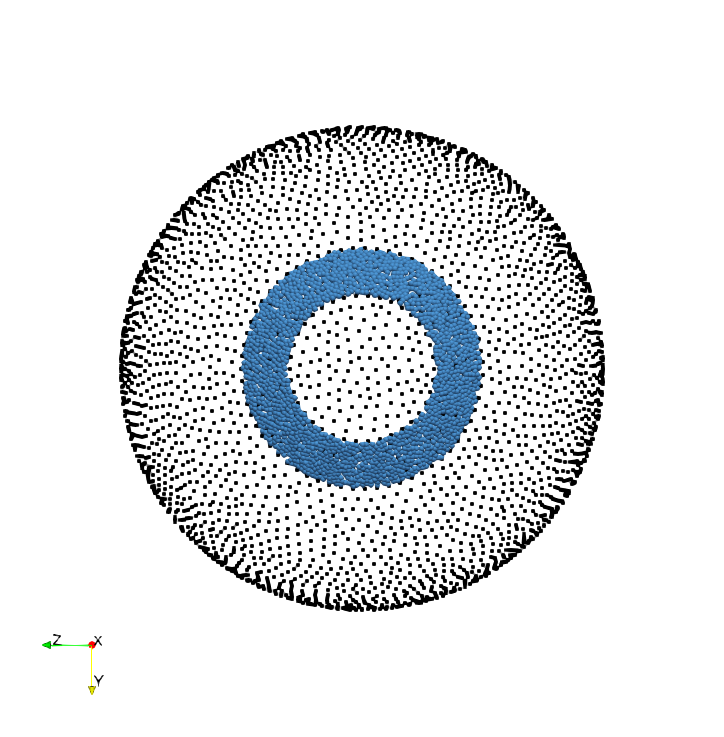}
  \caption{Domain and initial condition of the annular drop given by blue spherical glyphs, falling on a hemisphere marked with black points.}
  \label{Fig:DropFall_IC}%
\end{figure}

We now consider an example with free boundaries present. For this, we propose a benchmarking example similar to the drop falling standard case used in volumetric flow with free surfaces. The surface domain under consideration is given by a unit hemisphere with $x \leq 0$. An annular initial fluid droplet is considered as shown in Figure~\ref{Fig:DropFall_IC}. This is given by the region of intersection of the hemisphere and an annular cylinder given by $0.3^2 \leq y^2 + z^2 \leq 0.5^2$. Gravity is set to $0$, and the drop initiates movement with a prescribed initial velocity of $\vec{v}(\vec{x}, 0) = (y,-x,0)^T$.

\begin{figure}
  \centering
  \includegraphics[width=0.3\textwidth]{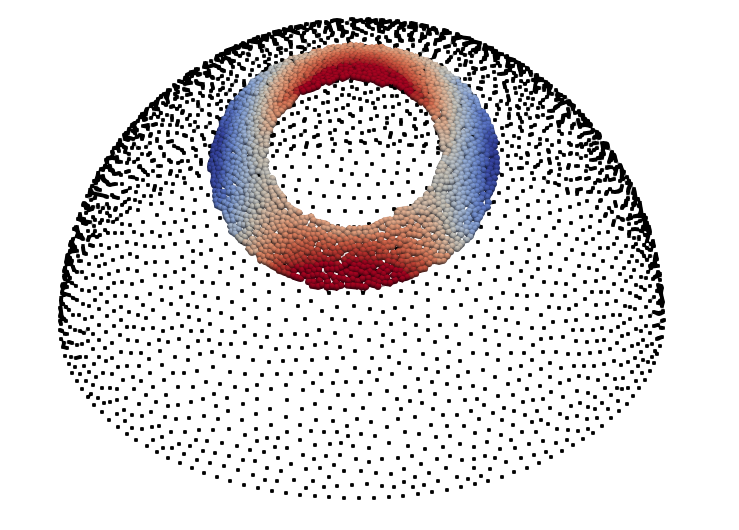}
  \includegraphics[width=0.3\textwidth]{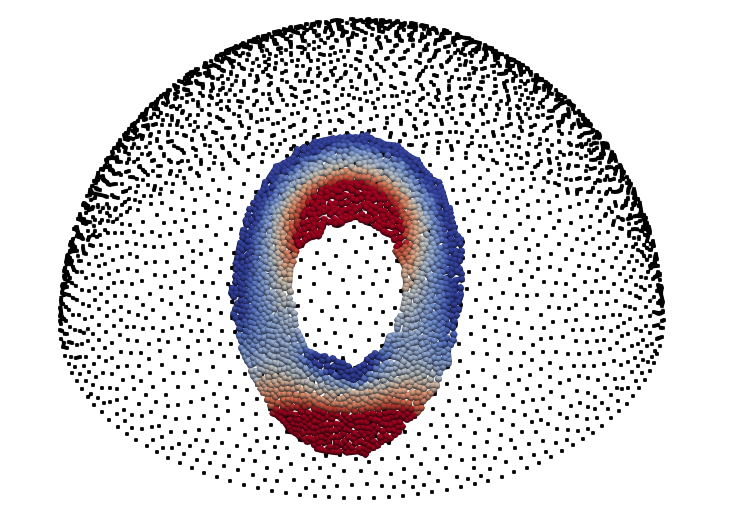}
  \includegraphics[width=0.1\textwidth]{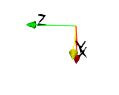}  
  \includegraphics[width=0.15\textwidth]{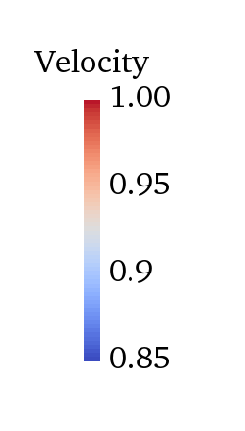}\\
  \includegraphics[width=0.3\textwidth]{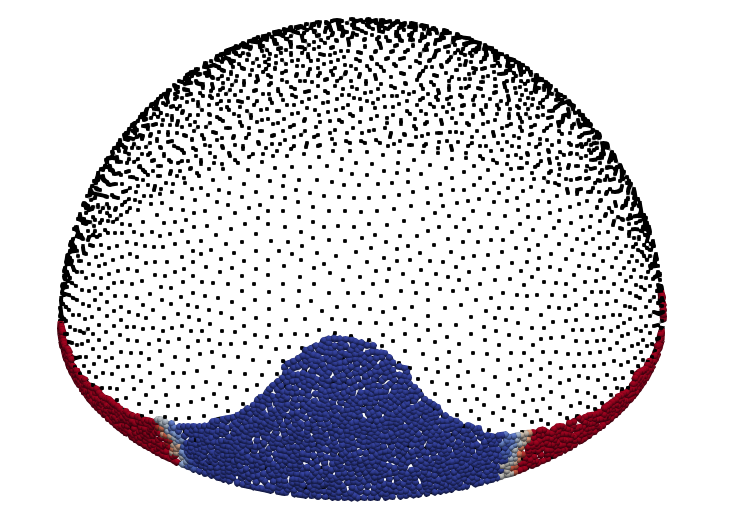}
  \includegraphics[width=0.3\textwidth]{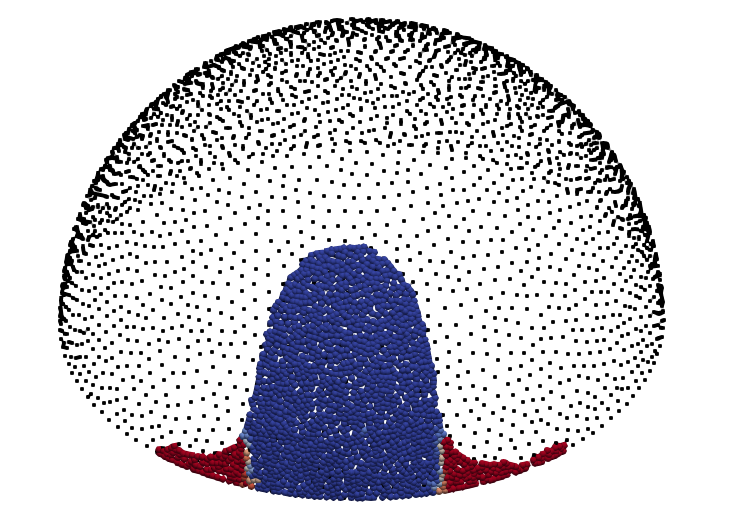}  
  \includegraphics[width=0.3\textwidth]{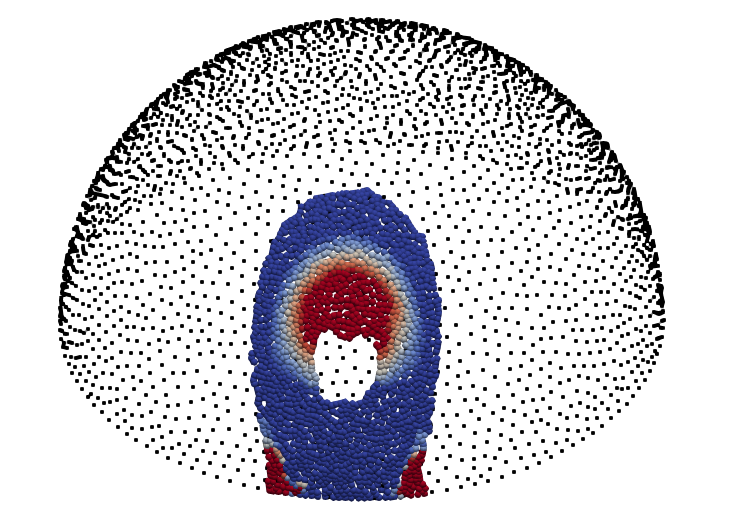}
  \caption{Falling annular droplet on a hemisphere. Domain at different times. Clockwise from top left: $t=0.33$~(top left), $t=0.86$~(top right), $t=1.13$~(bottom right), $t=1.4$~(bottom center), and  $t=1.88$~(bottom left). The fluid is shown with spherical coloured glyphs. And the surface is shown with black points. The colour of the fluid represents the velocity, coloured according to the initial velocity distribution.}
  \label{Fig:DropFall_Video}%
\end{figure}

The physical properties used are $\rho = 1$, $\eta = 0.001$. The initial point cloud is seeded with $h=0.06$, which corresponds $N=1607$ points. Furthermore, a fixed time step of $\Delta t = 0.006$ is used. Slip boundary conditions are enforced at the surface boundary. The evolution of the falling drop is shown in Figure~\ref{Fig:DropFall_Video}, with points coloured according to the initial velocity. 

\begin{figure}
  \centering
  \includegraphics[width=0.32\textwidth]{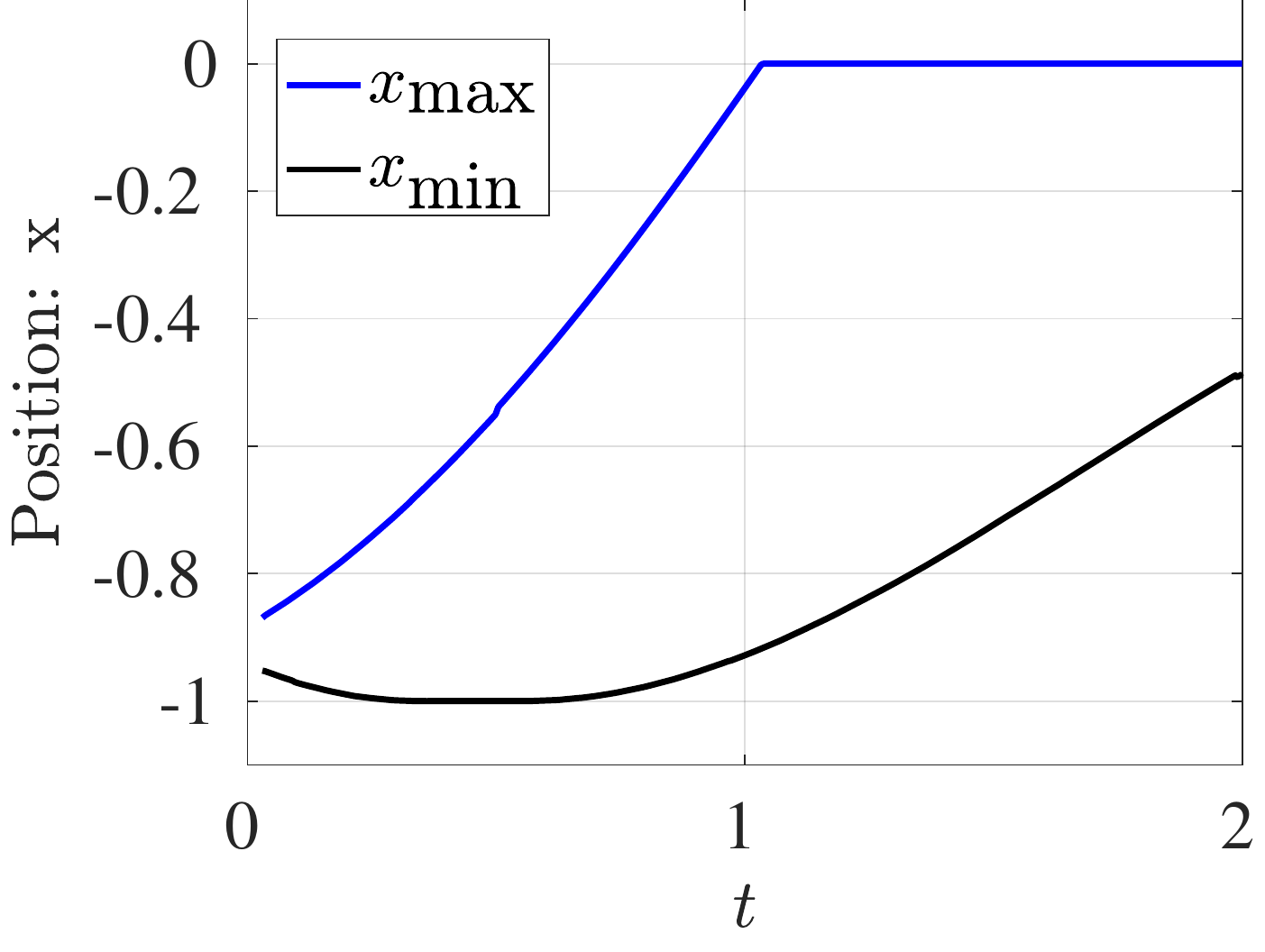}
  \includegraphics[width=0.32\textwidth]{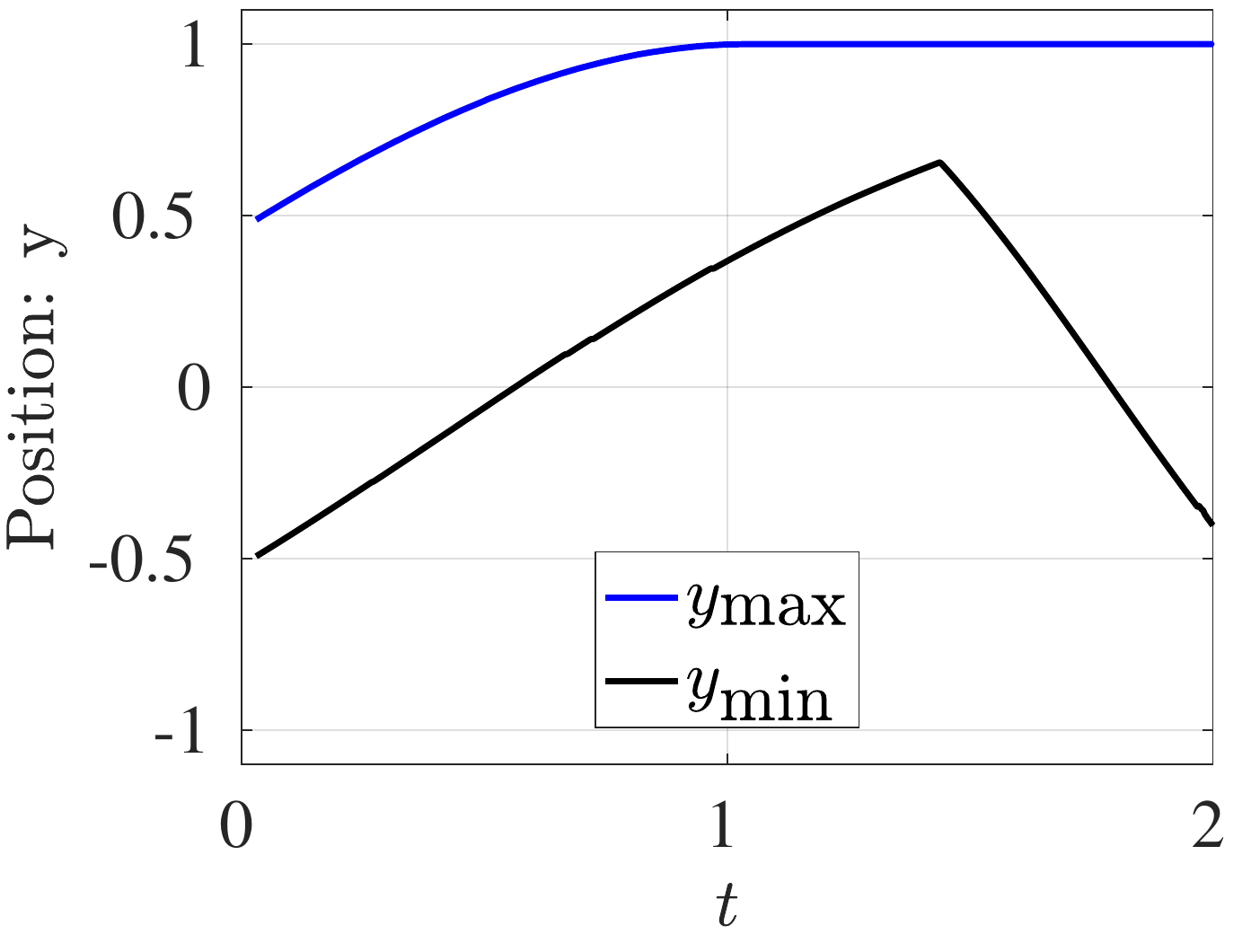}
  \includegraphics[width=0.32\textwidth]{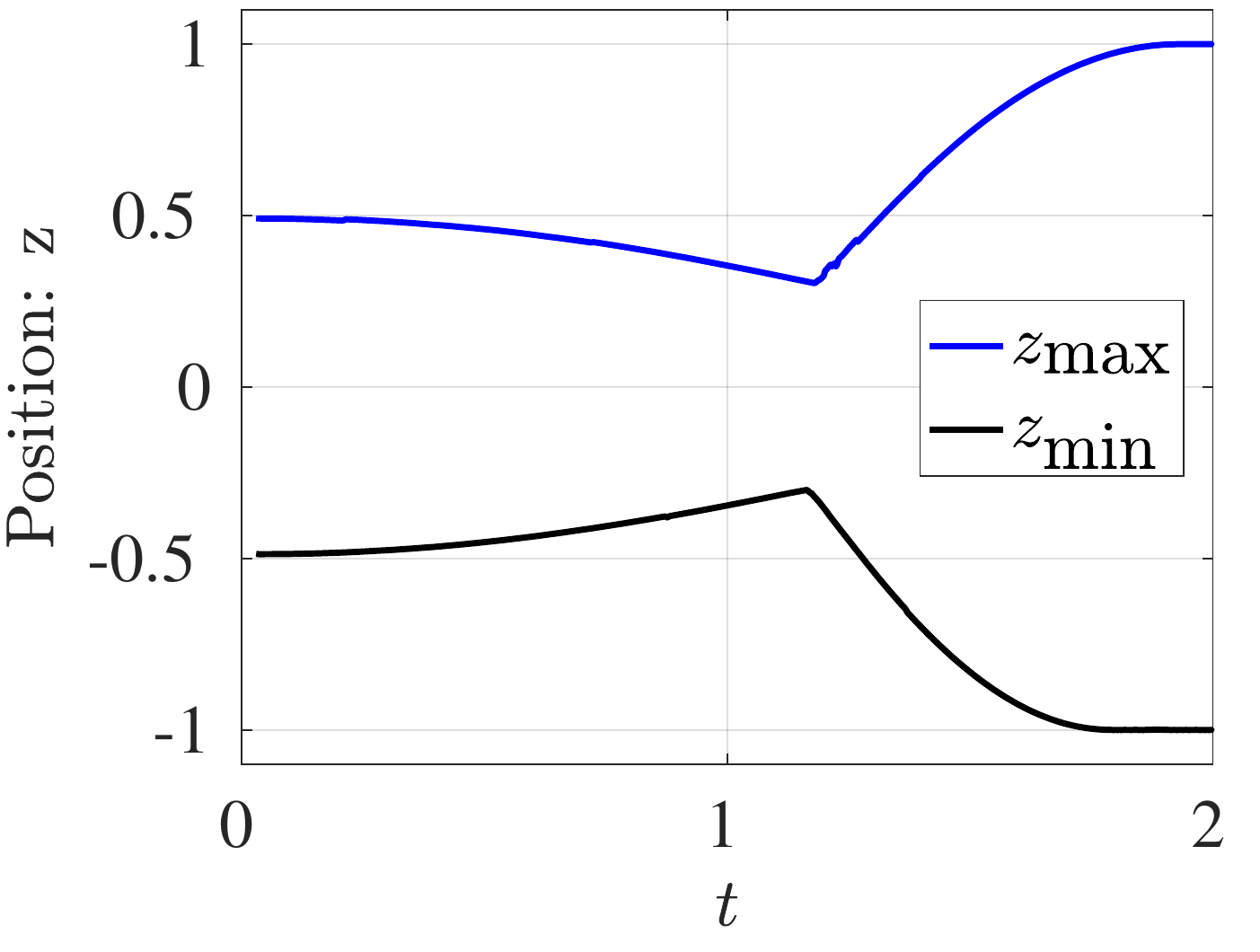}
  \caption{Drop falling: For the position vector $\vec{x}(t) = (x,y,z)^T$, the evolution of the maximum and minimum values of each coordinate are shown.}
  \label{Fig:DropFall_Position}%
\end{figure}

To quantify the results, we use similar metrics to those used in the classical dam break example used to benchmark free surface volumetric flows (e.g. \cite{Tiwari2003}). The evolution of the location of the droplet is quantified by plotting the maximum and minimum values of the $x$, $y$, and $z$ coordinates in Figure~\ref{Fig:DropFall_Position}. The prescribed initial velocity in the $y$ direction, with a $0$ velocity in the $z$ direction causes the drop to elongate in the $y$ direction. As a consequence of the mass conservation constraint, the drop thins along the $z$ direction as it falls. The annular drop reaches the boundary at $t=1.03$. The increase in the $z$ coordinate happens slightly after that. Just after hitting the boundary, the `hole' in the droplet starts becoming smaller, till it closes at $t=1.25$, after which the expansion in the $z$ direction starts. Further, the  maximum velocity reached is $\| \vec{v} \|_{\text{max}} = 2.16$.

This example also illustrates the free surface detections algorithm, and the importance of preventing faulty addition of points outside boundaries. Since no point is added at a location surrounded by existing boundary points, the hole in the droplet is not filled too early. 

It must be noted here that these results are not very representative physically as interface tension and contact angle effects have not been incorporated.

\subsection{Sloshing on an Armadillo}

\begin{figure}
  \centering
  \includegraphics[width=0.45\textwidth]{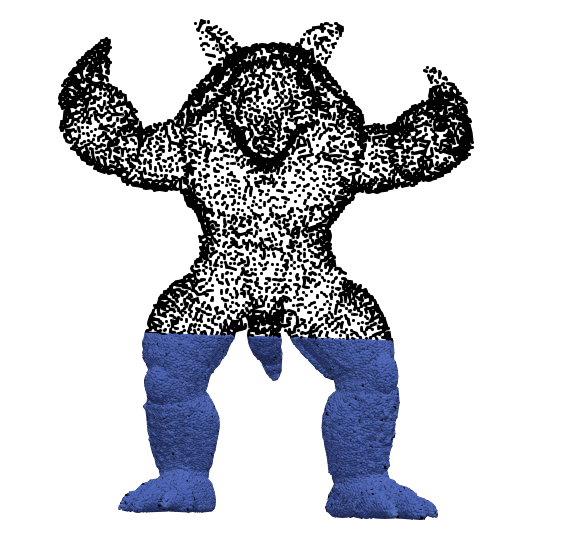}
  \includegraphics[width=0.45\textwidth]{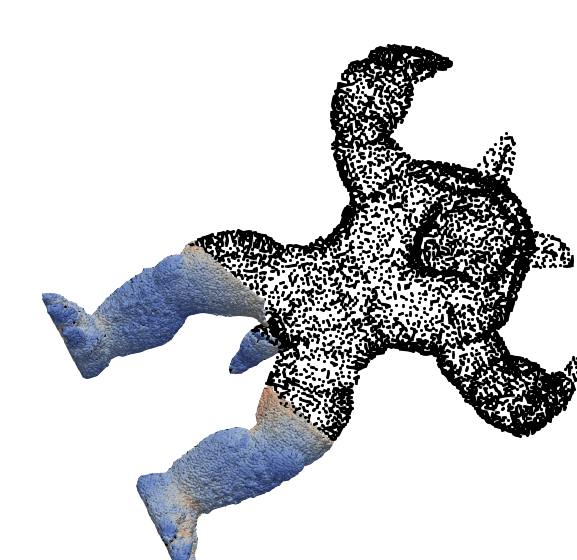}\\
  \includegraphics[width=0.45\textwidth]{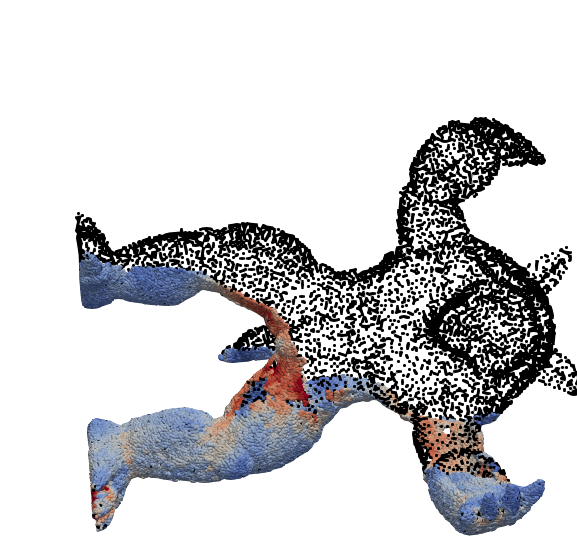}
  \includegraphics[width=0.45\textwidth]{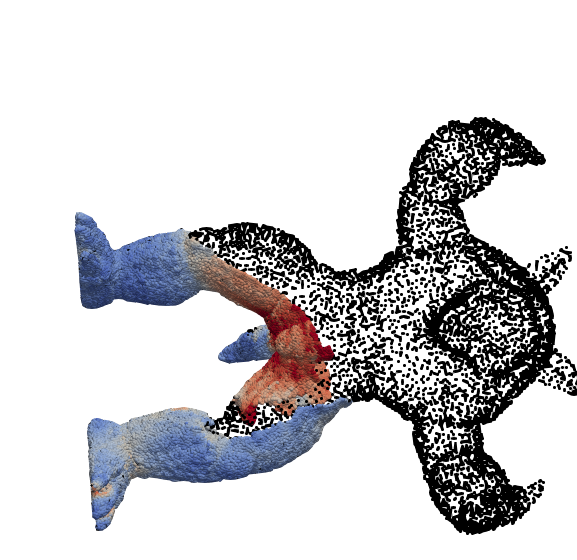}\\
  \includegraphics[width=0.45\textwidth]{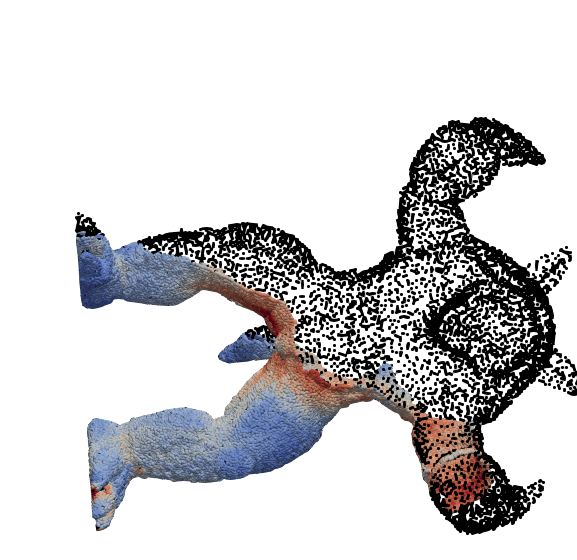}
  \includegraphics[width=0.45\textwidth]{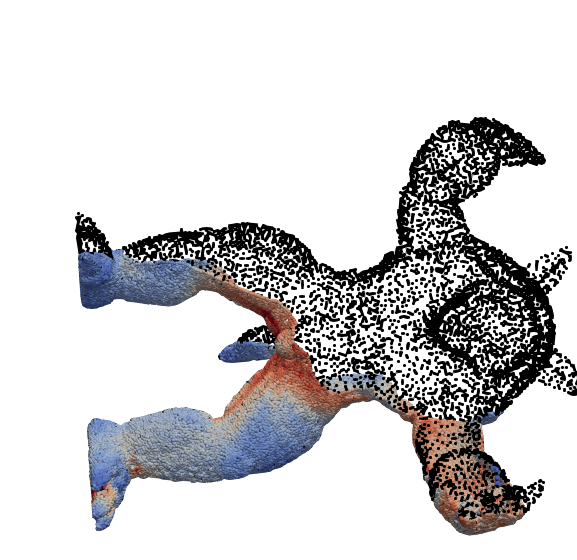}  
    \caption{Sloshing on an Armadillo surface. Clockwise from top left: $t=0$~(top left), $t=1.1$~(top right), $t=2.57$~(middle right), $t=5.31$~(middle left), $t=7.1$~(bottom left), and  $t=9.7$~(bottom right). The fluid is shown with spherical glyphs coloured by the velocity, and the surface is shown with black points.}
  \label{Fig:ArmadilloSlosh}%
\end{figure}

As the last example, we show a simulation with complex evolution of free boundaries where the fluid undergoes sloshing. The surface considered is the Stanford Armadillo \cite{Krishnamurthy1996}. At the initial time, fluid is filled on the armadillo surface up to a certain height, as shown in the first image of Figure~\ref{Fig:ArmadilloSlosh}. The armadillo is rotated till it becomes horizontal, which causes the fluid to fall due to gravity. 

A scaled version of the armadillo is used with a scale factor of $0.035$ to the original data set. The gravity field points downwards with respect to the initial state of the armadillo, and is prescribed by $\vec{g}=(0,1,0)^T$. $\rho = 1$ and $\eta = 0.01$ are used, with a varying time step of $\Delta t = 0.2 h / \|\vec{v}\|_{\infty}$. The evolution of the free boundaries are shown in Figure~\ref{Fig:ArmadilloSlosh}. The rotation of the armadillo surface causes the fluid to fall under gravity. Initially, fluid is present in three disconnected regions, which merge as the fluid falls. This setup illustrates the capability of this method to handle complex geometries and free boundary evolution. 
%


\section{Conclusion}
\label{sec:conclusion}

We introduced a meshfree Lagrangian framework to model fluid flow on curved and moving surfaces. The method is intrinsic to the surface in the sense that only the surface is discretized without any  form of discretization of the bulk around the surface. Meshfree points or particles are moved in a Lagrangian sense subject to the constraint of lying on the surface. The position of the surface itself can also be evolving in time. Methods to handle moving free boundaries within a surface were introduced, along with the boundary conditions for the same. Flow on a surface without free boundaries was investigated numerically using the basis of a point cloud that represents both the fluid domain, and the surface that constrained the fluid flow. On the other hand, for flow with free boundaries, the use of two point clouds was introduced. For the latter case, the main fluid point cloud was moved on another background point cloud that represented the (possibly moving) surface. The introduced Lagrangian framework incorporates local fixes for point cloud distortion, which removes the need for a possible global remeshing in moving mesh methods. As a result, the methods presented here provide a straightforward way to handle flow on moving and distorting surfaces, where the flow domain itself can also deform. 

Coupling this with a meshfree Generalized Finite Difference Method~(GFDM), we presented a strong form scheme to solve the Navier--Stokes equations on manifolds. This used a modified Chorin projection approach, where the tangential velocity constraint is enforced as an additional projection step, before the projection for incompressibility. Validation cases were introduced for the Lagrangian framework, and for incompressible surface Navier--Stokes flow with free boundaries. The numerical examples illustrate the potential of the introduced method to be used for a variety of applications of surface flow. 

Several directions of work could be pursued to improve the scheme presented here. Most importantly, different schemes need to be explored for improved stability with larger time steps, as has been done for meshfree collocation schemes for volumetric flow \cite{Suchde2018_INSE}. Particularly, the development of a scheme with a fully implicit velocity step, or an implicit coupled velocity-pressure scheme, need to be considered. The biggest challenge in this is handling the viscous term implicitly. Additionally, improvements in accuracy for the pressure field may also be obtained with an additional pressure Poisson equation at the end of the scheme, obtained by substituting the final velocity in the momentum equation, as has shown to be useful in meshfree GFDM schemes for volumetric flow \cite{Jefferies2015}. 

\end{document}